\documentclass[reqno]{amsart}
\usepackage[pdfusetitle]{hyperref}
\usepackage{preamble/packages}
\usepackage{preamble/colours}
\usepackage{preamble/mydiag}
\usepackage{preamble/environments}
\usepackage{preamble/notation}
\usepackage[top=30truemm,bottom=30truemm,left=35truemm,right=35truemm]{geometry}

\crefformat{section}{\S#2#1#3}
\crefmultiformat{section}{\S\S#2#1#3}{ and~#2#1#3}{, #2#1#3}{, and~#2#1#3}
\crefdefaultlabelformat{#2\textup{#1}#3}
\author{Samuel Lewis}
\title{Real variations of stability on K3 categories}
\date{}

\begin{document}

\begin{abstract}
  This paper proves that the \(2\)-Calabi--Yau triangulated category associated with the preprojective algebra \(\Ppa\) of an affine or hyperbolic graph admits many real variations of stability conditions, in the sense of Anno, Bezrukavnikov, and \text{Mirkovi\'{c}}. We do this by connecting the Coxeter arrangements of the graph with the homological algebra by introducing the concept of \emph{real flows}. This categorifies the notion of above and below for alcoves in a hyperplane arrangement and allows us to generalise much of the machinery of real variations for affine arrangements to the hyperbolic setting. In the process, we show that the derived category of \(\Ppa\)-modules with nilpotent cohomology is equivalent to the derived category of nilpotent \(\Ppa\)-modules.
\end{abstract}

  \maketitle
  \tableofcontents
    
  \section{Introduction}\label{sec:intro}

Let \(\Triang\) be a triangulated category and let \(\Hpln\) be a locally finite arrangement of real hyperplanes equipped with a direction \(\below\) on its alcoves (connected components of the complement). Introduced in \cite{anno_stability_2015}, a real variation of stability conditions on \(\Triang\), parametrised by \(\Hpln\) in the direction of \(\below\), consists of the following data:
\begin{enumerate}
  \item a real central charge \(\charge\),
  \item an assignment of hearts (of bounded t-structures on \(\Triang\)) to alcoves.
\end{enumerate}
This data needs to satisfy positive pairing and t-structure compatibility conditions for neighbouring alcoves, such that if an alcove is \say{above} another with respect to \(\below\), the associated hearts differ by an appropriate homological shift. Full details are given in \cref{dfn:realstab}. If the hearts are faithful, this compatibility can be thought of as a perverse equivalence of triangulated categories.

In this paper, we use the data of a finite connected graph \(\graph\) of affine or hyperbolic type, along with its preprojective algebra \(\Ppa\), to construct a large class of such stability conditions on the 2-Calabi--Yau triangulated category of \(\Ppa\)-modules with nilpotent cohomology
\begin{equation*}
  \Gcat*[\graph] \defn \Db[\Nilp]{\Ppa} \subset \Db{\Ppa},
\end{equation*}
parametrised by the Coxeter arrangement \(\Hpln[\graph]\) associated with \(\graph\). It is well-known that \(\Gcat*[\graph]\) and \(\Hpln[\graph]\) are closely related via tilting theory, with mutation functors categorifying simple reflections. As in the work \cite{ikeda_stability_2014} of Ikeda, we are motivated by the realisation of the space \(\Stab{\Gcat*[\graph]}\) of Bridgeland stability conditions on \(\Gcat*[\graph]\) as a covering space of some complex hyperplane arrangement coming from root data of \(\graph\). To avoid the technical issue of whether or not a functor that fixes all simples is the identity functor, we take an alternative approach to Bridgeland stability conditions via real variations.

\subsection{Main results}\label{ssec:result}

In the theory of quiver representations, much attention has been given to the finite and affine cases due to finite-dimensionality and noetherianity, respectively. This has important geometric consequences in dimension two, but loses meaning when the graph or quiver is wild. From the perspective of Coxeter combinatorics, it is known that in this last case a meaningful distinction can be made between the strictly wild case (indefinite) and the minimally wild case (hyperbolic). We are motivated to demonstrate the utility of the hyperbolic case in the aforementioned homological context, lifting affine stability conditions into a hyperbolic setting.

This is achieved via an explicit construction using the Coxeter arrangement \(\Hpln[\graph]\) associated with \(\graph\) and the tilting theory of \(\Ppa\). Our construction generalises the examples given in \cite{anno_stability_2015}. In \cref{dfn:realstab} we generalise the affine real variations of stability in \cite{anno_stability_2015} in two ways. Firstly, we allow for a central charge compatible with hyperbolic space, where we cannot easily translate hyperplanes to the \say{origin}; secondly, we categorify their notion of direction on alcoves away from an affine setting. 

In turn, this motivates us to introduce the new notion of a \emph{real flow} on a Coxeter arrangement, which we sketch here. Full details are given in \cref{ssec:flow}.

\begin{dfn}\label{dfn:realflow}
  Let \(\Hpln\) be a locally finite arrangement of real hyperplanes. A \emph{real flow} on \(\Hpln\) is a choice of \(\above\) or \(\below\) at every codimension one wall between neighbouring alcoves, such that around every codimension two wall, there is a unique \say{source} alcove \(\alc[\source]\), opposite which is a unique \say{target} alcove \(\alc[\target]\), as illustrated below.
  \begin{equation*}
    \begin{tikzpicture}[scale=3]
      \node[circle, fill, scale=.5] (F) at (0, 0) {};
      \node (As) at (0, 1/3) {\(\alc[\source]\)};
      \node (At) at (0, -1/3) {\(\alc[\target]\)};
      \node[right] (z1) at (1, 0) {};
      \node[above] (z2) at (1/2, {sqrt(3)/2}) {};
      \node[above] (z3) at (-1/2, {sqrt(3)/2}) {};
      \node[left] (z4) at (-1, 0) {};
      \node[below] (z5) at (-1/2, -{sqrt(3)/2}) {};
      \node[below] (z6) at (1/2, -{sqrt(3)/2}) {};
      \node[r, rotate=90] (p1) at ({2*sqrt(3)/9}, 0) {\Huge\(\prec\)};
      \node[r, rotate=150] (p2) at ({sqrt(3)/9}, 1/3) {\Huge\(\prec\)};
      \node[r, rotate=30] (p3) at (-{sqrt(3)/9}, 1/3) {\Huge\(\prec\)};
      \node[r, rotate=90] (p4) at (-{2*sqrt(3)/9}, 0) {\Huge\(\prec\)};
      \node[r, rotate=150] (p5) at (-{sqrt(3)/9}, -1/3) {\Huge\(\prec\)};
      \node[r, rotate=30] (p6) at ({sqrt(3)/9}, -1/3) {\Huge\(\prec\)};
      \tkzDefPoint(0,0){O}
      \tkzDefPoint(1,0){P}
      \tkzClipCircle(O,P)
      \hgline{0}{60}
      \hgline{60}{120}
      \hgline{120}{180}
      \hgline{180}{240}
      \hgline{240}{300}
      \hgline{300}{360}
      \hgline{0}{180}
      \hgline{60}{240}
      \hgline{120}{300}
    \end{tikzpicture}
    \begin{tikzpicture}[scale=3]
      \node[circle, fill, scale=.5] (F) at (0, 0) {};
      \node (As) at (0, 1/4) {\(\alc[\source]\)};
      \node (At) at (0, -1/4) {\(\alc[\target]\)};
      \node[above] (z1) at ({sqrt(2)/2}, {sqrt(2)/2}) {};
      \node[above] (z2) at (-{sqrt(2)/2}, {sqrt(2)/2}) {};
      \node[below] (z3) at (-{sqrt(2)/2}, -{sqrt(2)/2}) {};
      \node[below] (z4) at ({sqrt(2)/2}, -{sqrt(2)/2}) {};
      \node[r, rotate=135] (p1) at (1/5, 1/5) {\Huge\(\prec\)};
      \node[r, rotate=45] (p2) at (-1/5, 1/5) {\Huge\(\prec\)};
      \node[r, rotate=135] (p3) at (-1/5, -1/5) {\Huge\(\prec\)};
      \node[r, rotate=45] (p4) at (1/5, -1/5) {\Huge\(\prec\)};
      \tkzDefPoint(0,0){O}
      \tkzDefPoint(1,0){P}
      \tkzClipCircle(O,P)
      \hgline{45}{135}
      \hgline{135}{225}
      \hgline{225}{315}
      \hgline{315}{45}
      \hgline{45}{225}
      \hgline{135}{315}
    \end{tikzpicture}
  \end{equation*}
\end{dfn}

This is given a categorical interpretation in \cref{sec:real}, where the choice of above or below corresponds to a choice of positive or negative power of mutation functor between the hearts assigned to the alcoves. In fact, we demonstrate that choosing such a real variation amounts to choosing a real flow on \(\Hpln[\graph]\). Our first main result is that real flows on the Coxeter arrangement are very common.

\begin{ppn}[\ref{ppn:flow}]\label{ppn:realflow}
  Real flows \(\flow\) are precisely those assignments of mutation functors to wall-crossings that descend to the arrangement groupoid of \(\Hpln[\graph]\). Equivalently, \(\flow\) is a real flow if and only if braid relations hold around every codimension two wall.
\end{ppn}

Our generalised notion of real variations in \cref{dfn:realstab} removes all obstructions to lifting from the affine to the hyperbolic setting. After some work, \cref{ppn:realflow} leads to the following main result.

\begin{thm}[\ref{thm:rstab}]
  Let \(\graph\) be a connected graph of affine or hyperbolic type and choose a real flow \(\flow\) on its Coxeter arrangement \(\Hpln[\graph]\). Then there exists a real variation of stability conditions on \(\Gcat*[\graph]\), where the bounded t-structure at each alcove \(\alc\) is determined by the image of the corresponding Weyl group element under \(\flow\).
\end{thm}

This connects stability conditions with monodromy, in the sense that choosing a trivial monodromy homomorphism is the same as a real stability condition on \(\Gcat*[\graph]\).

Along the way, we prove a result that may be of independent interest regarding \emph{faithfulness} of the standard t-structure in \(\Gcat*[\graph]\), explored in the hereditary case in \cite[Lemma 3.2]{elagin_thick_2021}. Indeed, Bridgeland in \cite{bridgeland_stability_2009} raises the \say{slightly subtle} problem of whether \(\Gcat*[\graph]\) is equivalent to \(\Db{\Nilp[\Ppa]}\).  We answer this in the affirmative.

\begin{thm}[\ref{thm:faith}]
  Let \(\Ppa\) be the preprojective algebra of a graph without loops, which is not of finite type. Then there is an equivalence of categories
  \begin{equation*}
    \map*{\Db{\Nilp[\Ppa]}}{\Db[\Nilp]{\Ppa}}
  \end{equation*}
  induced by the inclusion \(\Nilp[\Ppa]\subset\Mod{\Ppa}\). That is, the standard heart \(\Nilp[\Ppa]\) is \emph{faithful} in \(\Gcat*[\graph]\).
\end{thm}

This allows us to rephrase the main compatibility condition of real variations in terms of perverse equivalences with respect to the filtration of \(\Gcat*[\graph]\) by thick triangulated subcategories generated by vertex simples.

\subsection{Structure of the paper}\label{ssec:struct}

In \cref{sec:back} we give the necessary preliminary material on graphs, Coxeter arrangements, and preprojective algebras, with an emphasis on hyperbolic type. In \cref{sec:real}, we intertwine the combinatorics and the homological algebra by defining the real central charge map \(\charge\) and t-structure assignment. The key notion in this section is that of a \emph{real flow}, a categorification of direction on \(\Hpln\). Following this, we are ready to prove our main result in \cref{sec:stab}, showing that the proposed construction satisfies the compatibility conditions for real variations.

\subsection{Notation and conventions}\label{ssec:conv}

All modules are right modules unless otherwise stated. We use the convention that paths in quivers are read from left to right. As is standard, we write \(\Db{\alg}\defn\Db{\Mod{\alg}}\) for a ring \(\alg\). The notation \(\hyp\) on a graph signifies an overextended ADE graph of the corresponding type, where the first subscript represents the number of vertices in the underlying ADE graph, and an additional subscript represents any extra edges added to this; all such notation is explained in \cref{ssec:data}. All subcategories are assumed to be strictly full, and the base field is \(\Comp\). We use the term \emph{flat} to mean a codimension two wall in the Coxeter arrangement, as we are not interested in flats of any other dimension. We denote by \(\thick[\simp[\ver]]\) the thick subcategory of \(\Gcat*[\graph]\) generated by \(\simp[\ver]\), so it is closed under shifts and exact triangles. Unless stated otherwise, graphs do not admit loops.

\subsection{Acknowledgements}\label{ssec:ack}

I wish to thank Franco Rota for sharing his expertise in stability conditions, Matthew Pressland for many helpful discussions on quivers and preprojective algebras, and my supervisors Gwyn Bellamy and Michael Wemyss for their constant support throughout this project.

\subsection*{Funding} The author was partly supported by EPSRC, and partly by ERC Consolidator Grant 101001227 (MMiMMa). 

\subsection*{Open access} For the purpose of open access, the author has applied a Creative Commons Attribution (CC:BY) licence to any Author Accepted Manuscript version arising from this submission.

  \section{Background}\label{sec:back}

\subsection{Input and output}\label{ssec:data}

Let \(\graph\) be a finite connected graph with vertex set \(\Ver{\graph}\), edge set \(\Edge{\graph}\), and denote its \(\emph{rank}\) by \(\rank\defn\card{\Ver{\graph}}\). To \(\graph\) we can associate a symmetric \(\rank\times\rank\) integer matrix
\begin{equation*}
  \GCM[\graph] = \paren{\GCM[\ver\ver*]} \defn 2\idmat - \adj[\graph]
\end{equation*}
known as its \emph{generalised Cartan matrix} (GCM), where \(\idmat\) is the identity matrix and \(\adj\) is the adjacency matrix of \(\graph\). Usually we assume that \(\graph\) has no loops, but if \(\graph\) has \(\val[\ver]\) loops at vertex \(\ver\) then we may take \(\GCM[\ver\ver] = 2-2\val[\ver]\).

More generally, GCMs correspond to Dynkin diagrams, and most of the ideas in this section still apply to GCMs that are not symmetric but are \emph{symmetrisable}. Recall that the definiteness of \(\GCM[\graph]\) splits \(\graph\) into three types:
\begin{enumerate}
  \item finite (spherical) type --- \(\GCM[\graph]\) is positive definite,
  \item tame (affine) type --- \(\GCM[\graph]\) is positive semidefinite,
  \item wild (indefinite) type --- \(\GCM[\graph]\) is indefinite.
\end{enumerate}
In the context of graphs without orientation (symmetric GCMs), finite type graphs are ADE and affine type graphs are extended ADE. Of special interest are graphs that are \say{barely} wild, as these provide a sensible first step into non-noetherian behaviour.

\begin{dfn}\label{dfn:hbolic}
  A wild graph \(\graph\) is said to be of \emph{hyperbolic type} if every proper full subgraph \(\graph'\subset\graph\) is of finite or affine type.
\end{dfn}

The name \emph{hyperbolic} comes from the signature of \(\GCM[\graph]\), which necessarily has \(\rank-1\) positive eigenvalues and one negative eigenvalue whenever \(\graph\) is hyperbolic. However, this is not a sufficient condition for \(\graph\) to be hyperbolic; in general, more is needed.

Hyperbolic Dynkin diagrams, of which \cref{dfn:hbolic} is a special case, were first classified in the works of \cite{li_classification_1988} and \cite{saclioglu_dynkin_1989}, and they have a rank of at most ten. A more recent classification for \(3\leq\rank\leq 10\), which fixes some errata, is contained in \cite{carbone_classification_2010}. As might be expected, a hyperbolic graph can be constructed by taking an affine ADE diagram and adding an \say{overextended vertex} to the extended vertex, but there are also exceptional cases.

\begin{ppn}\label{ppn:hbolic}
  If \(\graph\) is a graph of hyperbolic type, then it is either one of the overextended ADE diagrams
  \begin{equation*}
    \begin{tikzcd}[every arrow/.append style={dash}, column sep = .9em, row sep = .9em]
      &[-0em] \dote \ar[d] \ar[ddl] \ar[ddr] && \Inf \ar[r] & \dote \ar[dr] &&&& \dote \ar[dl] &&& \Inf \ar[r] & \dote \ar[d] \\
      & \Inf && && \dote \ar[dl] \ar[r] & \cdots\cdots \ar[r] & \dote \ar[dr] &&&&& \dote \ar[d] \\
      \dote \ar[r] & \cdots\cdots \ar[r] &[-0em] \dote & & \dote &&&& \dote && \dote \ar[r] & \dote \ar[r] & \dote \ar[r] & \dote \ar[r] & \dote
    \end{tikzcd}
  \end{equation*}
  \begin{equation*}
    \begin{tikzcd}[column sep = 2em]
      \hypA[1,\dots,7] &&&&& \hypD[4,\dots,8] &&&&& \hypE[6]
    \end{tikzcd}
  \end{equation*}
  \begin{equation*}
    \begin{tikzcd}[every arrow/.append style={dash}, column sep = .8em, row sep = .8em]
      &&&& \dote \ar[d] &&&&&&& \dote \ar[d] \\
      \Inf \ar[r] & \dote \ar[r] & \dote \ar[r] & \dote \ar[r] & \dote \ar[r] & \dote \ar[r] & \dote \ar[r] & \dote && \dote \ar[r] & \dote \ar[r] & \dote \ar[r] & \dote \ar[r] & \dote \ar[r] & \dote \ar[r] & \dote \ar[r] & \dote \ar[r] & \Inf
    \end{tikzcd}
  \end{equation*}
  \begin{equation*}
    \begin{tikzcd}[column sep = 5em]
        \hypE[7] &&& \hypE[8] = \finE[10]
    \end{tikzcd}
  \end{equation*}
  or a member of one of the following infinite families
  \begin{equation*}
		\begin{tikzcd}[every arrow/.append style={dash}, column sep = 1.2em]
			\dote 
			\ar[loop, in = -30, out = 30, distance = 2em, "1", swap]
			\ar[loop, in = -30, out = 30, distance = 5em, "n\geq2"]
			& \cdots &&& \dote \ar[rr, draw=none, "\raisebox{1ex}{\vdots}" description]
			\ar[rr, bend left, "1"]
			\ar[rr, bend right, swap, "n\geq3"]
			&& \dote
			&& \dote \ar[rr, draw=none, "\raisebox{1ex}{\vdots}" description]
			\ar[rr, bend left,        "1"]
			\ar[rr, bend right, swap, "n\geq1"]
			\ar[loop, out = 138, in = -135, distance = 2em]
			&& \dote
			&& \dote \ar[rr, draw=none, "\raisebox{1ex}{\vdots}" description]
			\ar[rr, bend left,        "1"]
			\ar[rr, bend right, swap, "n\geq1"]
			\ar[loop, out = 138, in = -135, distance = 2em]
			&& \dote \ar[loop, out = -48, in = 45, distance = 2em] \\
			& \Loop[n] &&&& \Kronecker[n] &&&& \Kronecker[n, 1] &&&& \Kronecker[n ,2]
		\end{tikzcd}
	\end{equation*}
  or one of the following exceptional graphs
  \begin{equation*}
		\begin{tikzcd}[every arrow/.append style={dash}, column sep = 1.3em, row sep = 1em]
      &&&&&[-0em] \dote \ar[ddl, shift left =.5] \ar[ddl, shift right =.5] \ar[ddr] &&&& \dote \ar[ddl, shift left =.5] \ar[ddl, shift right =.5] \ar[ddr] &&&& \dote \ar[ddl, shift left =.5] \ar[ddl, shift right =.5] \ar[ddr, shift left =.5] \ar[ddr, shift right =.5] \\[-0.5*0.57735em]
      \\[-0.5*0.28868em]
			\Inf \ar[r, shift left =.5] \ar[r, shift right =.5] & \dote \ar[r, shift left =.5] \ar[r, shift right =.5] & \dote && \dote \ar[rr] && \dote && \dote \ar[rr, shift left =.5] \ar[rr, shift right =.5] && \dote && \dote \ar[rr, shift left =.5] \ar[rr, shift right =.5] && \dote
		\end{tikzcd}
	\end{equation*}
  \begin{equation*}
    \begin{tikzcd}
      \hypA[1,1] &&& \affA[2,1] &&&[-1.5em] \affA[2,2] &&& \affA[2,3] 
    \end{tikzcd}
  \end{equation*}
  \begin{equation*}
    \begin{tikzcd}[every arrow/.append style={dash}, column sep = 1.3em, row sep = 1em]
      &[-0em] \dote \ar[d] \ar[ddl] \ar[ddr] &&&& \dote \ar[ddl] \ar[ddr] \ar[d] &&&& \dote \ar[ddl] \ar[ddr] \ar[d] &&&& \dote \ar[dd] \ar[ddl] \ar[ddr] \ar[dl] \ar[dr] \\[-0.5*0.57735em]
      & \Inf \ar[dl] &&&& \Inf \ar[dl] \ar[dr] &&&& \Inf \ar[d] &&& \dote && \dote \\[-0.5*0.28868em]
      \dote \ar[rr] &&[-0em] \dote && \dote \ar[rr] && \dote && \dote \ar[r] & \dote \ar[r] & \dote && \dote & \dote & \dote
    \end{tikzcd}
  \end{equation*}
  \begin{equation*}
    \begin{tikzcd}
      \hypA[2,1] &&& \hypA[2,2] &&&[-1.5em] \hypA[3,1] &&& \Star[5] 
    \end{tikzcd}
  \end{equation*}
  where the \(\infty\) nodes correspond to overextended vertices in the sense that the full subgraph without this node is extended ADE.
\end{ppn}
\begin{proof}
  The stated diagrams with \(\rank\geq 3\) can all be found in \cite[Tables 1--23]{carbone_classification_2010}. They are the diagrams with the \(\rank\times\rank\) identity matrix in the penultimate column, meaning that their symmetrising matrix is trivial, and we identify symmetric double arrows with a pair of edges. Consider \(\rank=1\). Trivially, any proper subgraph of \(\Loop[\val]\) is of finite type, and for \(\val\geq 2\), this is a wild graph. For \(\rank = 2\), we can have at most one loop at each vertex as \(\Loop[1]=\affA[0]\) is tame and \(\Kronecker[\val, 1(2)]\) are wild for all \(\val\geq 1\).
\end{proof}

For future reference, the graphs in \cref{ppn:hbolic} are summarised, by rank, in the following table.
\begin{equation}\label{tbl:hbolic}
	\begin{array}{c | c}
		\rank & \graph \\
		\hline
		1 & \Loop[\val] \\
		2 & \Kronecker[\val] \quad \Kronecker[\val,1] \quad \Kronecker[\val,2] \\
		3 & \hypA[1] \quad \hypA[1,1] \quad \affA[2,1] \quad \affA[2,2] \quad \affA[2,3] \\
		4 & \hypA[2] \quad \hypA[2,1] \quad \hypA[2,2] \\
		5 & \hypA[3] \quad \hypA[3,1] \\
		6 & \hypA[4] \quad \hypD[4] \quad \Star[5] \\
		7 & \hypA[5] \quad \hypD[5] \\
		8 & \hypA[6] \quad \hypD[6] \quad \hypE[6] \\
		9 & \hypA[7] \quad \hypD[7] \quad \hypE[7] \\
		10 & \hypD[8] \quad \hypE[8]
	\end{array}
\end{equation}

Future sections study stability conditions arising from the extended (\say{affine}) ADE graphs, as well as those in \eqref{tbl:hbolic}. Perhaps unsurprisingly, Coxeter arrangements constructed in \cref{ssec:cox} sit in affine \(\rank\)-space for the former and hyperbolic \(\rank\)-space for the latter. In \cref{ssec:cox}, the necessity of considering only these graphs and no wilder is made clear.

\begin{dfn}[Cf.\ {\cite[Definition 1]{anno_stability_2015}}]\label{dfn:realstab}
  Let \(\Triang\) be an (Ext-finite) triangulated category and let \(\rspace\) be a finite-dimensional Euclidean or hyperbolic space containing a discrete collection of hyperplanes \(\Hpln\). Let \(\Alc\) be the set of \emph{alcoves}, the connected components of \(\rspace\minus\Hpln\). Suppose also that \(\Alc\) carries the notion of \say{above/below} for neighbouring alcoves \(\alc,\alc*\), denoted \(\below\). A \emph{real variation of stability conditions} on \(\Triang\), parametrised by \(\rspace\minus\Hpln\) in the direction of \(\below\), consists of the following:
  \begin{enumerate}[label=(\alph*)]
    \item A real (analytic) map \(\map[\charge]{\rspace}{\Kgrp{\Triang}[\Real]\dual},\ \stab\mapsto\charge[\stab]\), the \emph{central charge},
    \item a map \(\map[\hrt]{\Alc}{\set{\Abel\subseteq\Triang\text{ bounded heart}}}\),
  \end{enumerate}
  satisfying the following properties for all \(\alc\in\Alc\).
  \begin{enumerate}
    \item If \(\mod\in\Abel\defn\hrt(\alc)\) is nonzero, then \(\pair{\charge[\stab]}{\Kthy{\mod}}>0\) for all \(\stab\in\alc\).
    \item If \(\alc*\above\alc\) is an alcove neighbouring \(\alc\) across a hyperplane \(\hpln\in\Hpln\), consider for \(\deg\in\Int[1]\) the full subcategories
    \begin{align*}
      \Abel[\hpln]\filt[\deg]&\defn\set{M\in\Abel}[\pair{\charge[\blank]}{\Kthy{M}}\colon \rspace\to\Real \text{ has a zero of order at least } \deg \text{ along } \hpln] \subseteq \Abel, \\
      \Gcat*[\alc,\hpln]\filt*[\deg] &\defn\set{X\in\Triang}[\coH[\Abel]{\sdot}{X}\in\Abel[\hpln]\filt[\deg]] \subseteq\Triang.
    \end{align*}
    Then \(\Abel*\defn\hrt(\alc*)\) is compatible with the filtration \(\Gcat*[\alc,\hpln]\filt*\) of \(\Triang\) in the sense that
    \begin{equation}\label{eqn:intprop}
      \Abel*[\hpln]\filt[\deg] = \Abel*\inter\Gcat*[\alc,\hpln]\filt*[\deg]
    \end{equation}
    \begin{equation}\label{eqn:quoprop}
      \frac{\Abel*[\hpln]\filt[\deg]}{\Abel*[\hpln]\filt[\deg+1]} = \frac{\Abel[\hpln]\filt[\deg]}{\Abel[\hpln]\filt[\deg+1]}\shift{\deg} \inside* \quo{\Gcat}{\Gcat*[\alc,\hpln]\filt*[\deg]}
    \end{equation}
    for all \(\deg\in\Int[0]\), where \(\Abel*[\hpln]\filt[\deg]\) is defined analogously to \(\Abel[\hpln]\filt[\deg]\).
  \end{enumerate}
\end{dfn}

The cohomology above is taken with respect to \(\Abel\), that is, \(\coH[\Abel]{\deg}{\obj}\defn\trunc[\Abel]{0}\trunc*[\Abel]{0}\obj\shift{\deg}\), where the truncation functors are adjoints to the inclusion of the aisle and coaisle of the bounded t-structure with heart \(\Abel\).

Notice that the first axiom gives a map into the space of Bridgeland stability conditions by taking a \say{real slice} of the usual central charge, and the second axiom connects neighbouring alcoves as \(\rspace\minus\Hpln\) is disconnected. Ours is a slight generalisation of the original definition, which has a specific affine notion of direction on the arrangement and a polynomial, rather than analytic, central charge. The definition in \cite{anno_stability_2015} assigns bounded t-structures themselves to the alcoves, but as bounded t-structures are determined by their hearts, our assignment is equivalent. Notably, the original definition works over a real vector space, which can be thought of as making a choice of origin in \(\rspace\). We lack a way to make this choice in hyperbolic space, motivating the categorical approach in \cref{sec:real}. 

\subsection{Coxeter framework}\label{ssec:cox}

Given \(\graph\) as above, consider the \emph{root lattice} \(\Rts\) and its real dual \(\Ths\), namely
\begin{equation*}
  \Rts\defn\direct*[\ver\in\Ver{\graph}]\Int\rt[\ver], \qquad \Ths \defn\direct*[\ver\in\Ver{\graph}]\Real\rt[\ver]\dual,
\end{equation*}
so that \(\Ths = \Hom[\Int]{\Rts}{\Int}\tensor[\Int]\Real\). The \(\rt[\ver]\) are called \emph{simple roots} and the \(\rt[\ver]\dual\) form the dual basis. Then \(\Rts\) admits a symmetric bilinear form given by
\begin{equation*}
  \form{\rt[\ver]}{\rt[\ver*]}\defn\GCM[\ver\ver*], \qquad \ver,\ver*\in\Ver{\graph}.
\end{equation*}

Recall that the \emph{Weyl group} of \(\graph\) is the Coxeter group \(\Weyl[\graph]\), with generating set
\begin{equation*}
  \set{\sref[\ver]}[\ver\in\Ver{\graph} \text{ without loops}]
\end{equation*}
and relations \(\sref[\ver]<2> = 1\), \(\paren{\sref[\ver]\sref[\ver*]}<2> = 1\) if \(\GCM[\ver\ver*]=0\), \(\paren{\sref[\ver]\sref[\ver*]}<3> = 1\) if \(\GCM[\ver\ver*]=-1\), and \(\paren{\sref[\ver]\sref[\ver*]}<\infty> = 1\) if \(\GCM[\ver\ver*]\leq-2\). There is an associated \emph{braid group} \(\Br[\graph]\), which has the same generating set and relations, except that the generators are no longer involutions. The Weyl group is finite if and only if \(\graph\) is of finite type, see \cite{kac_infinite-dimensional_1990}. There is an action of \(\Weyl\) on \(\Rts\) via the \emph{simple reflections}
\begin{equation*}\label{eqn:srfl}
  \sref[\ver]\vect \defn \vect - \form{\vect}{\rt[\ver]}\rt[\ver],
\end{equation*}
giving rise to a \emph{root system} in \(\Rts\). The \emph{real roots} are the image of the simple roots under this action, and for each such \(\rt\) there is a corresponding hyperplane
\begin{equation*}
  \hpln[\rt] \defn \set{\stab\in\Ths}[\pair{\stab}{\rt} = 0], \qquad \hpln[\ver] \defn \hpln[\rt[\ver]]
\end{equation*}
in \(\Ths\). The group \(\Weyl\) acts on \(\Ths\) by linearly extending the dual action. The standard (Weyl) chamber is the upper orthant
\begin{equation*}
  \Ths[+] \defn \set{\stab\in\Ths}[\pair{\stab}{\rt[\ver]} > 0 \text{ for all } \ver\in\Ver{\graph}].
\end{equation*}

\begin{dfn}\label{dfn:tits}
  The \emph{Tits cone} associated with \(\graph\) is the \(\rank\)-dimensional convex cone
  \begin{equation*}
    \TC\defn\union*[\weyl\in\Weyl]\weyl\cl{\Ths[+]}\subseteq\Ths.
  \end{equation*}
\end{dfn}

The region \(\TC\) is the whole of \(\Real<\rank>\) if and only if \(\graph\) is finite type (\cite[Proposition 3.12(e)]{kac_infinite-dimensional_1990}), as there are no imaginary roots. Regardless, the classical result \cite[Theorem 1]{moody_root_1979} implies that if \(\graph\) is of finite, affine, or hyperbolic type, then its real and imaginary roots are given by
\begin{equation*}
  \set{\vect\in\Rts}[\form{\vect}{\vect} = 2], \qquad \set{\vect\in\Rts}[\form{\vect}{\vect} \leq 0],
\end{equation*}
respectively. Thus, in the \say{at worst hyperbolic} setting, we can dualise and obtain
\begin{equation}\label{eqn:titscone}
  \TC\union-\TC = \set{\stab\in\Ths}[\stab\tran(\Adj\GCM[\graph])\stab \geq 0].
\end{equation}
This explicit description of the Tits cone enables us to take a codimension one \say{slice} of the Tits cone to more easily visualise the combinatorics and see connections to spherical, affine, and hyperbolic geometry. In the theory of affine root systems, this role is played by the \emph{affine level}
\begin{equation*}
  \Level[\aff] \union-\Level[\aff] \defn \set{\stab\in\Ths}[\pair{\stab}{\mir}<2> = 1] \subset \TC,
\end{equation*}
where \(\mir\in\Rts\) is the \emph{minimal imaginary root}. We generalise this to the hyperbolic setting as follows. For affine type, the nullspace of \(\GCM[\graph]\) is given as
\begin{equation}\label{eqn:gcmnull}
  \set{\vect\in\Rts}[\vect\tran\GCM[\graph]\vect = 0] = \Int\set{\mir},
\end{equation}
see \cite[\S 2.6]{humphreys_reflection_1990}. Then \(\GCM[\graph]\Adj\GCM[\graph] = \det\GCM[\graph]\idmat[\rank] = 0\), and so \((\Adj\GCM[\graph])\stab\) is in the nullspace of \(\GCM[\graph]\) for all \(\stab\in\Ths\). Using \eqref{eqn:gcmnull}, we can conclude that
\begin{equation*}\label{eqn:formaff}
  \stab\tran(\Adj\GCM[\graph])\stab = \uplambda\pair{\stab}{\mir}<2>.
\end{equation*}
The constant of proportionality \(\uplambda\) is the (positive) determinant of the corresponding finite type Cartan matrix, which is \(\rank+1\) in type \(\finA\), \(4\) in type \(\finD\), and \(9-\rank\) in type \(\finE\). In this way, the following level set of \eqref{eqn:titscone} naturally generalises \(\Level[\aff]\).

\begin{dfn}\label{dfn:level}
  The \emph{real level} associated with a graph \(\graph\) of affine or hyperbolic type is the \((\rank-1)\)-dimensional region  
  \begin{equation*}
    \Level[\graph] \defn \set*{\stab\in\Ths}[\sum\nolimits_{\ver\in\Ver{\graph}}\stab[\ver]\geq 0,\ \stab\tran(\Adj\GCM[\graph])\stab = 1].
  \end{equation*}
\end{dfn}

Intersecting each real root hyperplane with the real level gives the hyperplane arrangement
\begin{equation*}
  \Hpln[\graph] \defn \union*[\rt\in\rsys[+]](\hpln[\rt]\inter\Level[\graph]),
\end{equation*}
which we call the \emph{Coxeter arrangement} of \(\graph\). This is a locally finite (every point of \(\Level\) is contained in at most finitely many hyperplanes) and essential (minimal intersections of
hyperplanes are points) arrangement of hyperplanes. The \emph{alcoves} of the Coxeter arrangement are the connected components of the complement \(\Level\setminus\Hpln\), denoted \(\csa\reg\) in \cite{anno_stability_2015}.

\begin{ntn}\label{ntn:alc}
  Write \(\Alc\) for the set of alcoves, which are in one-to-one correspondence (see, for example, \cite[Propositions 8.23 and 8.27]{hall2013lie}) with Weyl group elements as an image of the \emph{standard alcove} \(\fAlc \defn \Level\inter\Ths[+]\).
\end{ntn}

If \(\graph\) is of hyperbolic type, then we may equivalently use the inverse of \(\GCM[\graph]\) rather than its adjugate, flipping the sign and obtaining the defining equation of a hyperboloid in \(\rank-1\) dimensions. This can then be identified with the hyperbolic space \(\Hyperbolic[\rank-1]\), which is most easily visualised in the disc model. For more detail on this see, for example, \cite{cannon_hyperbolic_1997}.

\begin{exm}[\(\rank = 2\)]\label{exm:rktwo}
  Suppose that \(\graph=\Kronecker[\val]\) for some \(\val\in\Int[1]\), namely
  \begin{equation*}
    \graph = \begin{tikzcd}[every arrow/.append style={dash}, column sep = 1.2em]
      0 \ar[rr, draw=none, "\raisebox{1.25ex}{\vdots}" description]
			\ar[rr, bend left, "1"]
			\ar[rr, bend right, swap, "\val"]
			&& 1,
    \end{tikzcd} \qquad
    \GCM[\graph] = \Matrix*{2&-\val\\-\val&2}, \qquad \Adj[\GCM[\graph]] = \Matrix*{2&\val\\\val&2}, \qquad \det\GCM[\graph] = 4-\val<2>,
  \end{equation*}
  where the GCM has eigenvalues \(2+\val\) and \(2-\val\). Note that the case \(\val=1\) is \(\finA[2]\), whilst \(\val=2\) is \(\affA[1]\). Now, \(\Level\union-\Level\) is the set of \(\stab\in\Ths\iso\Real<2>\) satisfying
  \begin{equation*}
    \stab[1]<2> + \val\stab[1]\stab[2] + \stab[2]<2> = 1.
  \end{equation*}
  This is an ellipse for \(\val = 1\), a pair of straight lines for \(\val = 2\), and a hyperbola for \(\val\geq 3\). In the latter two cases, this is asymptotically parallel to the Tits cone, which has slopes
  \begin{equation*}
    \Kcon[\val]<\pm1> \defn \frac{\val\pm\sqrt{\val<2>-4}}{2},
  \end{equation*}
  a value that we call the \emph{Kronecker constant} and satisfies \(\Kcon[\val]+\Kcon[\val]\inv = \val\). Note that \(\Kcon[2] = 1\), \(\Kcon[3]\) is the square of the golden ratio, and \(\Kcon[\val]\) tends to \(\infty\) so that \(\TC\) limits to \(\Ths[+]\). The Weyl groups are isomorphic for \(\val=2\) and \(\val\geq 3\) (the free group on the generators \(\sref[0], \sref[1]\)), so the Coxeter arrangements are the same up to the size of the alcoves.
\end{exm}

\begin{figure}[ht]
  \centering 
  \includegraphics[width=.33\textwidth]{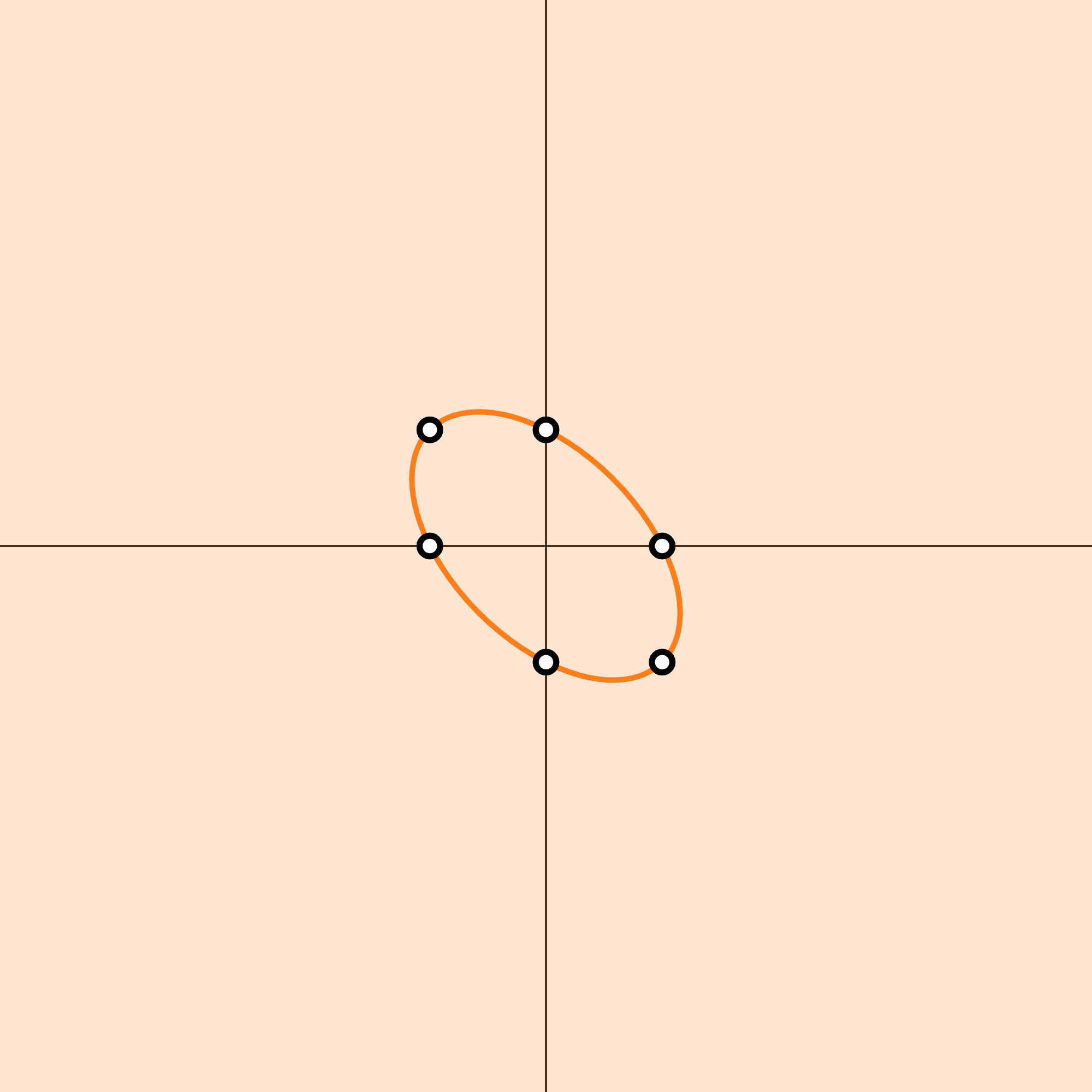}\hfill
  \includegraphics[width=.33\textwidth]{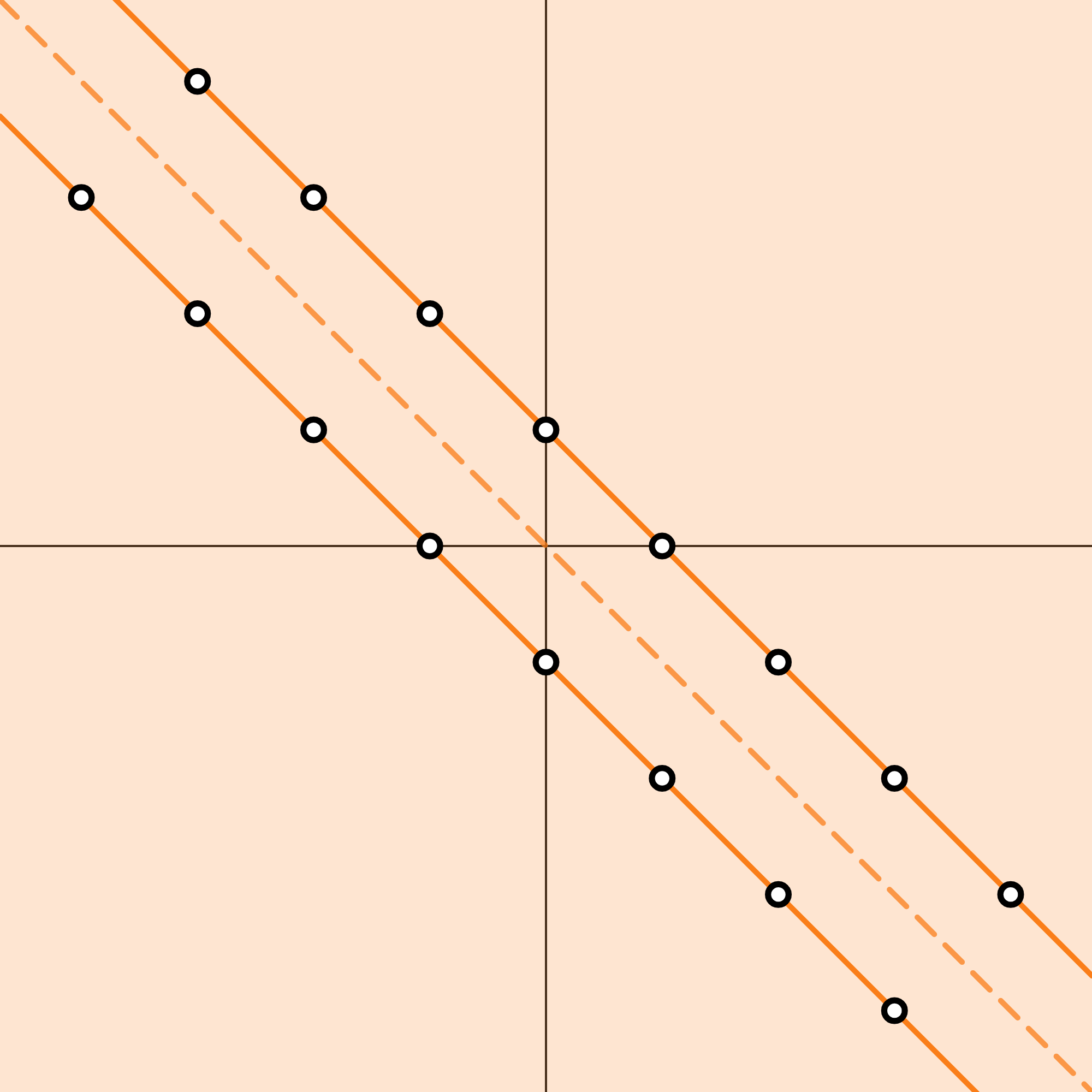}\hfill
  \includegraphics[width=.33\textwidth]{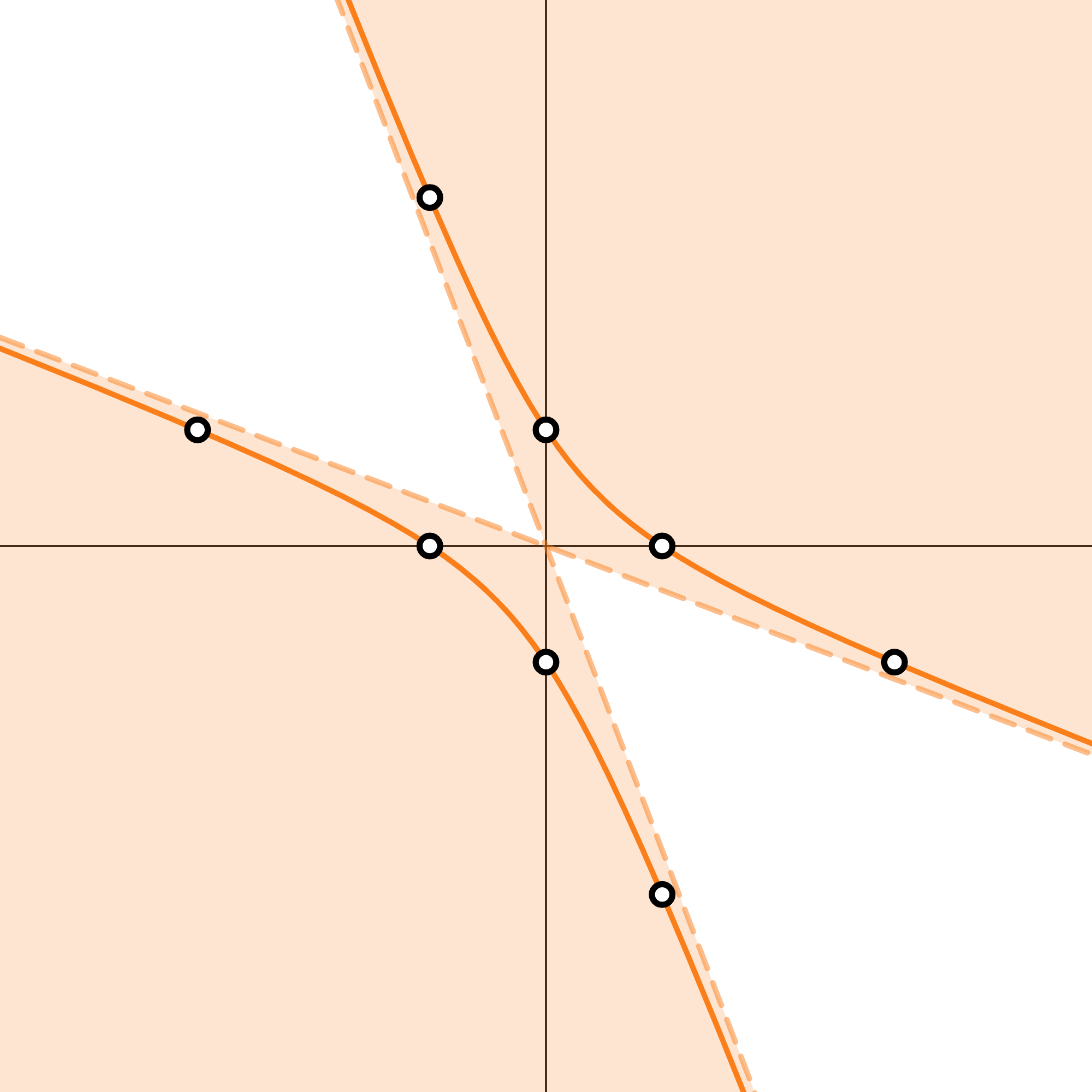}
  \caption{\(\TC\union-\TC\) (shaded region with dashed boundary) and \(\Level\union-\Level\) (solid line) for \(\Kronecker[1]\), \(\Kronecker[2]\), and \(\Kronecker[3]\), showing the transition between the spherical, affine, and hyperbolic setting. Open dots represent intersections with real root hyperplanes \(\hpln[\rt]\), of which there are finitely many for \(\Kronecker[1]\).}
\end{figure}

\begin{exm}[\(\rank = 3\), affine]\label{exm:rkthreeaff}
  Suppose that \(\graph = \affA[2]\), namely 
  \begin{equation*}
    \graph = \begin{tikzcd}[every arrow/.append style={dash}, column sep = 1em, row sep = 1em]
      & 0 \ar[dl] \ar[dr] \\
      1 \ar[rr] && 2
    \end{tikzcd}, \qquad
    \GCM[\graph] = \Matrix*{2&-1&-1\\-1&2&-1\\-1&-1&2}, \qquad \Adj[\GCM[\graph]] = \Matrix*{3&3&3\\3&3&3\\3&3&3}, \qquad \det\GCM[\graph] = 0.
  \end{equation*}
  Note that \(\det\GCM[\finA[2]] = 3\), as in \eqref{eqn:formaff} with \(\mir=\paren{1,1,1}\). Thus, \(\Level\union-\Level\) is the set of \(\stab\in\Ths\) satisfying
  \begin{equation*}
    \paren{\stab[0] + \stab[1] + \stab[2]}<2> = 1
  \end{equation*}
  in \(\Ths\), and as usual, we take \(\Level[\graph]\) to be the points living above \(\stab[0] + \stab[1] + \stab[2] = 0\). Intersecting this with the root hyperplanes gives the standard honeycomb tiling of \(\Real<2>\). The Weyl group is the affine symmetric group \(\Sym*[3]\).
\end{exm}

\begin{exm}[\(\rank = 3\), hyperbolic]\label{exm:rkthreehyp}
  Suppose that \(\graph = \hypA[1]\), namely
  \begin{equation*}
    \graph = \begin{tikzcd}[every arrow/.append style={dash}, column sep = 1.3em]
      \infty \ar[r] & 0 \ar[r, shift left =.5] \ar[r, shift right =.5] & 1,
    \end{tikzcd} \qquad
    \GCM[\graph] = \Matrix*{2&-1&0\\-1&2&-2\\0&-2&2}, \qquad \Adj[\GCM[\graph]] = \Matrix*{0&2&2\\2&4&4\\2&4&3}, \qquad \det\GCM[\graph] = -2.
  \end{equation*}
  Then the real level is the following hyperboloid in \(\Real<3>\)
  \begin{equation*}
    \Level[\graph] = \set{\stab\in\Ths}[\stab[\infty]+\stab[0]+\stab[1]\geq0,\ 2\stab[0]<2> + \tfrac{3}{2}\stab[1]<2> + 2\stab[\infty]\stab[0] + 2\stab[\infty]\stab[1] + 4\stab[0]\stab[1] = 1].
  \end{equation*}
  By a standard change of coordinates using the eigenvalues \(2\), \(2 + \sqrt{5}\), and \(2-\sqrt{5}\) of \(\GCM[\graph]\), we may diagonalise this to obtain the standard hyperboloid
  \begin{equation*}
    -\stab*[\infty]<2> + \stab*[0]<2> + \stab*[1]<2> = -1,
  \end{equation*}
  which can be used to give an explicit parametric form for the real central charge defined in \cref{ssec:charge}. The Weyl group is generated by the involutions \(\sref[\infty], \sref[0], \sref[1]\), subject to the relations
  \begin{equation*}
    \sref[\infty]\sref[0]\sref[\infty] = \sref[0]\sref[\infty]\sref[0], \qquad \sref[\infty]\sref[1] = \sref[1]\sref[\infty].
  \end{equation*}
  There is no relation between \(\sref[0]\) and \(\sref[1]\). This can be seen in the Coxeter arrangement of \(\hypA[1]\), which contains \(4\)-gons, \(6\)-gons, and \(\infty\)-gons in its dual graph.
\end{exm}

When \(\graph\) is hyperbolic, \(\Level\) can be identified with the hyperbolic space \(\Hyperbolic[\rank-1]\), which is most easily seen via the hyperboloid model. For more details on this, see \cite{cannon_hyperbolic_1997}. Proceeding in the same way as \cref{exm:rkthreehyp}, the remaining rank three hyperbolic graphs \(\hypA[1,1]\), \(\affA[2,1]\), \(\affA[2,2]\), and \(\affA[2,3]\) give rise to four more tilings of the hyperbolic plane. All rank three Coxeter arrangements coming from graphs of affine or hyperbolic type are shown below.

\begin{figure}[ht]
  \begin{equation*}
    \begin{array}{ccc}
      \includegraphics[width=0.275\textwidth]{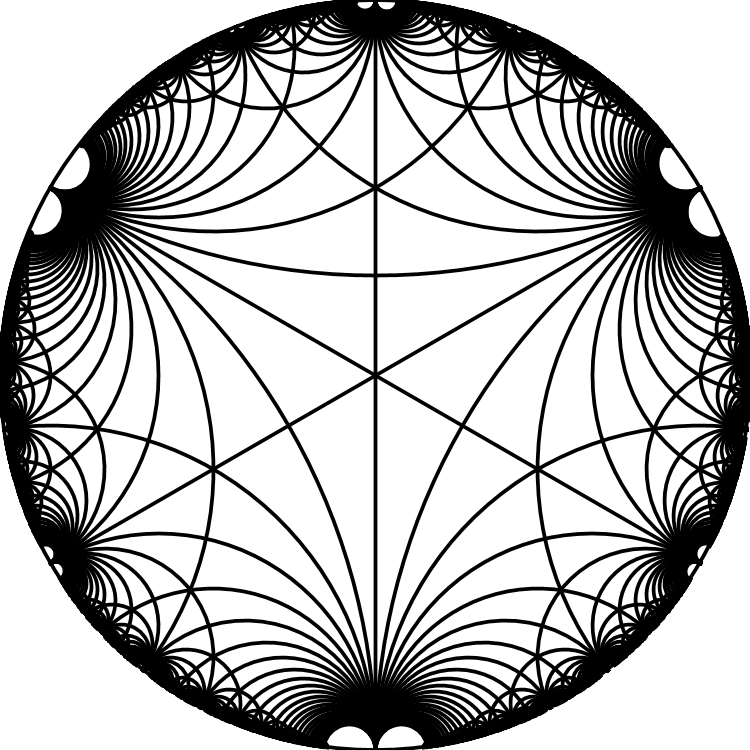} & \includegraphics[width=0.324\textwidth]{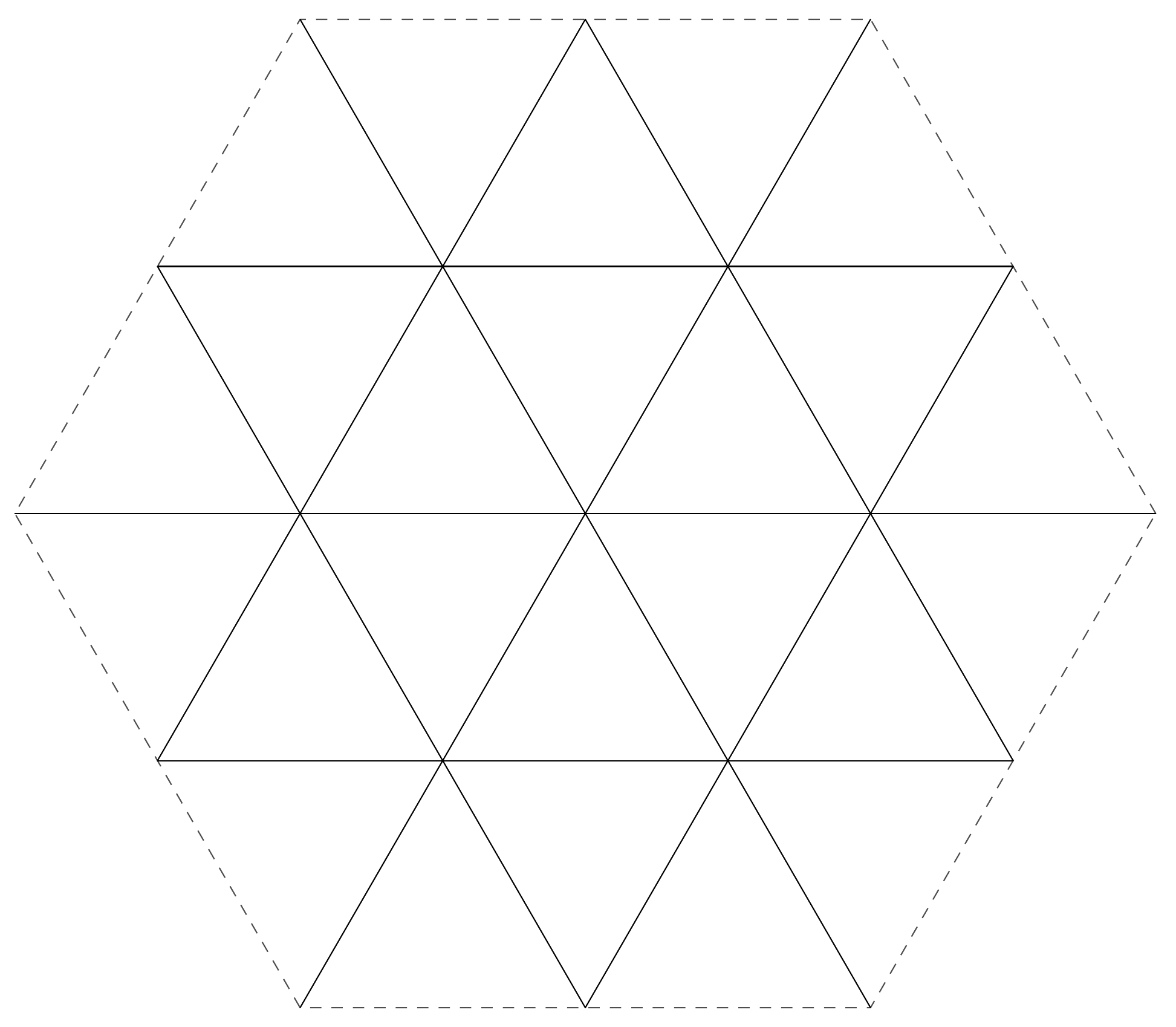} & \includegraphics[width=0.275\textwidth]{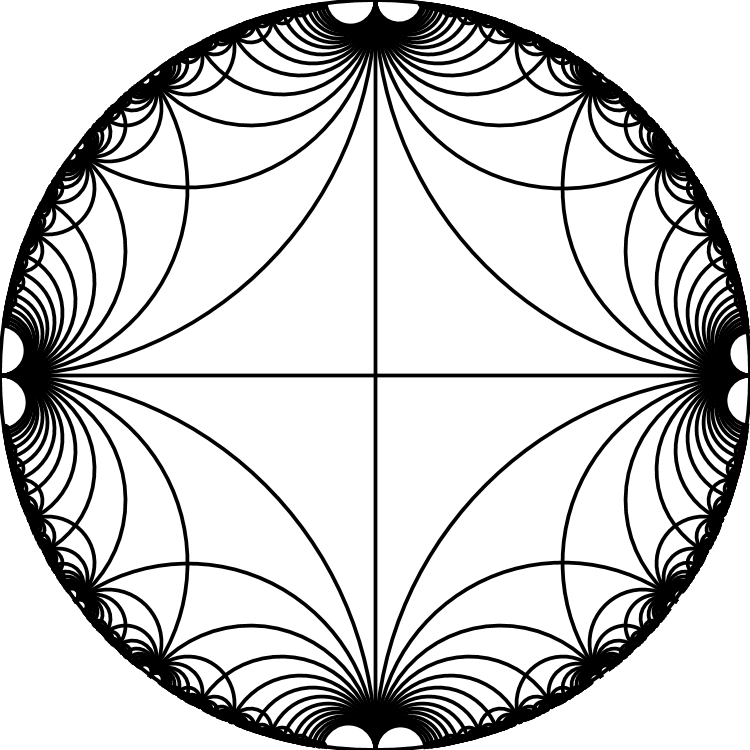} \\
      \hypA[1]\ (2, 3, \infty) & \affA[2]\ (3, 3, 3) & \hypA[1,1]\ (2, \infty, \infty) \\
      \includegraphics[width=0.275\textwidth]{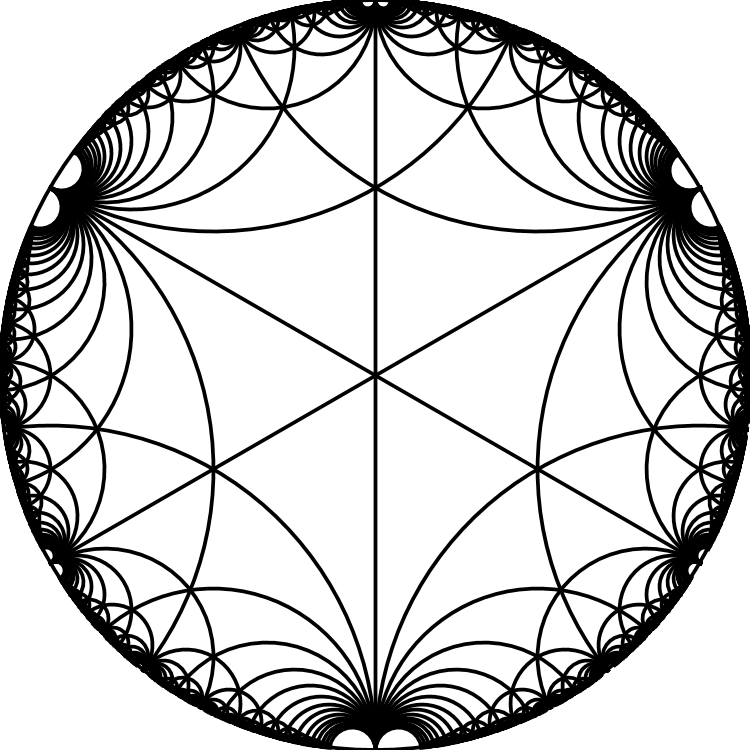} & \includegraphics[width=0.275\textwidth]{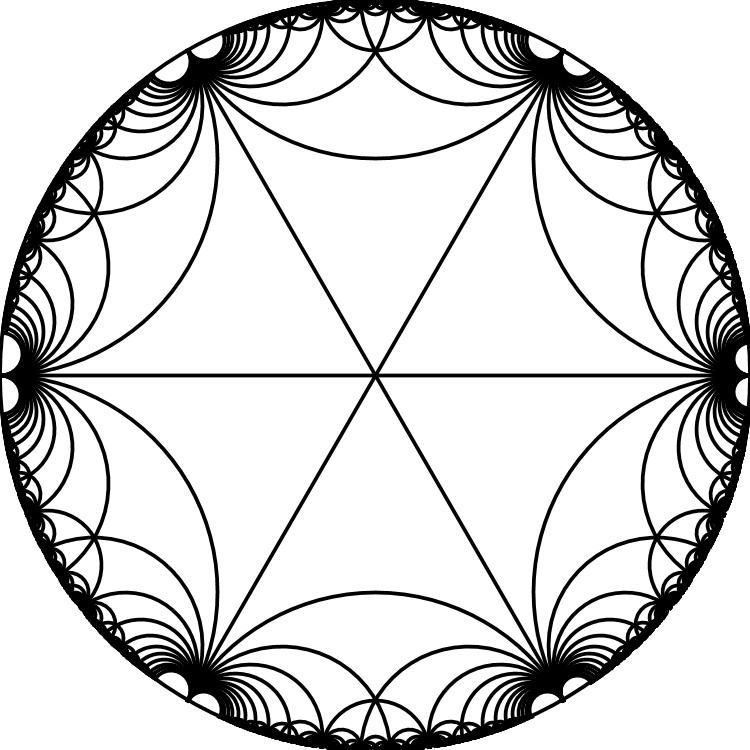} & \includegraphics[width=0.275\textwidth]{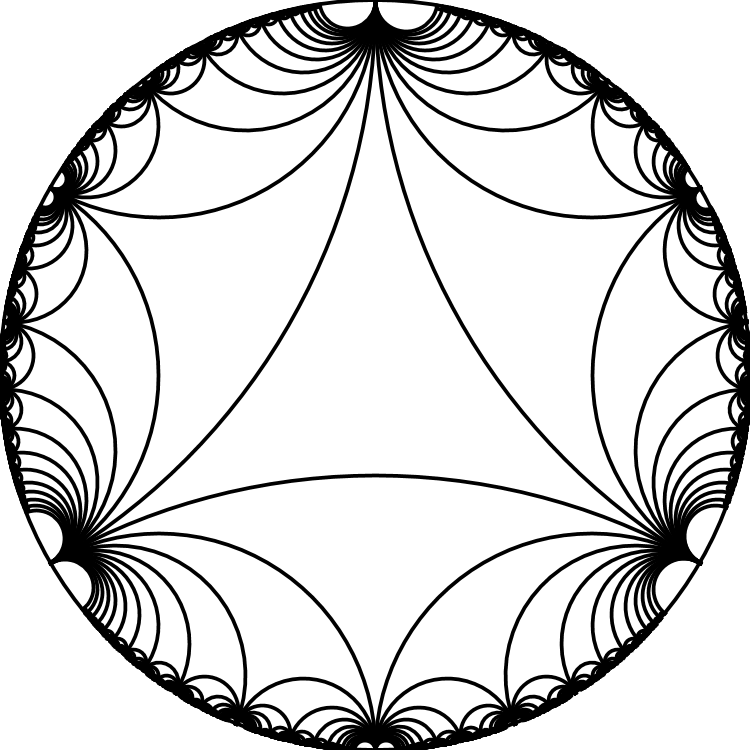} \\
      \affA[2,1]\ (3, 3, \infty) & \affA[2,2]\ (3, \infty, \infty) & \affA[2,3]\ (\infty, \infty, \infty)
    \end{array}
  \end{equation*}
  \caption{Rank three Coxeter arrangements of affine or hyperbolic type.}
\end{figure}

\subsection{Arrangement groupoids}\label{ssec:gdel}

This section recalls the various categories that can be associated with the Coxeter arrangement \(\Hpln[\graph]\) for any affine or hyperbolic \(\graph\). We follow \cite[\S 1.4]{deligne_les_1972}. The Coxeter arrangement \(\Hpln[\graph]\) admits a categorical description via the \emph{free category} \(\Free\Hpln[\graph]\), which has an object for each alcove and is freely generated by \say{wall-crossing} morphism pairs \((\map{\alc}{\alc*}, \map{\alc*}{\alc})\) for each pair of adjacent alcoves \(\alc,\alc*\in\Alc\). This can be thought of as a quiver path algebra on the dual arrangement to \(\Hpln[\graph]\) and hence also goes by the name \emph{Deligne quiver}. A concatenation of morphisms \(\upbeta\) is called a \emph{positive path}, and we say that such a path is \emph{reduced} if it doesn't cross the same hyperplane twice. Since \(\Hpln[\graph]\) is locally finite, this is the same as \(\map[\upbeta]{\alc[1]}{\alc[2]}\) having minimal length among all paths between the alcoves \(\alc[1]\) and \(\alc[2]\); see \cite[Lemma 2]{salvetti_topology_1987}. 

Now consider the smallest equivalence relation \(\sim\) on \(\Free\Hpln[\graph]\) that is compatible with morphism composition, that identifies reduced positive paths \(\map[\path,\path*]{\alc}{\alc*}\) with the same source and target. Then the \emph{path category} is defined as
\begin{equation*}
  \Path\Hpln[\graph] \defn \quo{\Free\Hpln[\graph]}{\!\sim},
\end{equation*}
often called the category of positive paths. If \(\alc\) and \(\alc*\) are two alcoves in a Coxeter arrangement, shown in a rank four (dimension three) hyperbolic cartoon below,
\begin{equation*}
  \begin{tikzpicture}
    \node[inner sep=0pt] (A) at (0,0) {$\alc$\includegraphics[scale=.2]{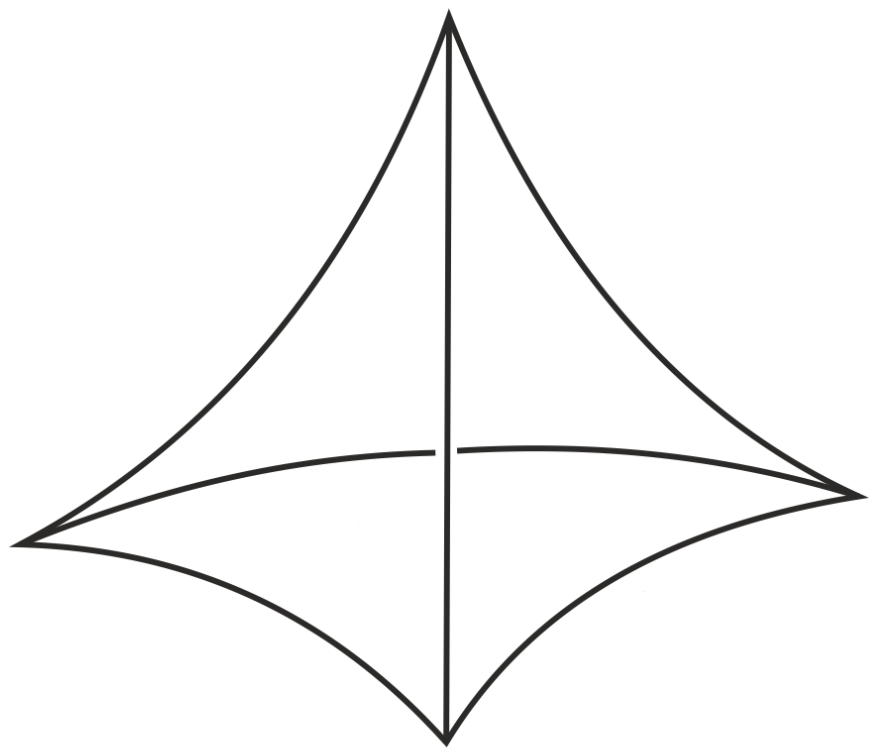}};
    \node[inner sep=0pt] (B) at (6,2) {\scalebox{-1}[1]{\includegraphics[scale=.15]{image/alcove.png}}$\alc*$};
    \draw [->] (A) edge[bend left] node[midway, above left] {\(\weyl = \sref[\ver[\len]]\cdots\sref[\ver[1]]\)} (B) ;
    \draw [->] (A) edge[bend right] node[midway, below right] {\(\weyl = \sref[\ver*[\len]]\cdots\sref[\ver*[1]]\)} (B) ;
  \end{tikzpicture}
\end{equation*}
then as \(\Weyl\) acts simply transitively on the set of alcoves, there is a unique \(\weyl\in\Weyl[\graph]\) such that \(\weyl\alc = \alc*\). However, in general, there are multiple reduced expressions for this \(\weyl\), which the following theorem connects. 

\begin{thm}[{\cite{matsumoto_generateurs_1964}}]\label{thm:matsu}
  Let \(\weyl[1], \weyl[2]\) be two reduced expressions for some \(\weyl\in \Weyl\). Then \(\weyl[1], \weyl[2]\) differ only by a (finite) sequence of braid relations.
\end{thm}

In other words, if \(\weyl[1],\weyl[2]\in\Br\) are reduced and identified in \(\Weyl\), then \(\weyl[1] = \weyl[2]\) in the braid group. 

\begin{dfn}[{\cite[(1.25)]{deligne_les_1972}}]\label{dfn:gdel}
  The \emph{arrangement groupoid} of \(\Hpln[\graph]\) is the groupoid completion \(\gDel[\graph]\) of \(\Path\Hpln[\graph]\). This simply means that every morphism in \(\Path\Hpln[\graph]\) is formally inverted.
\end{dfn}

\subsection{Homological framework}\label{ssec:hom}

Given a graph \(\graph\), we may consider its \emph{doubled quiver}, namely an oriented graph \(\double{\graph}\) with vertex set \(\Ver{\graph}\) but with two arrows \(\edge, \edge\rev\) (with the opposite orientation) for each edge in \(\Edge{\graph}\). Write \(\source(\edge)\) and \(\target(\edge)\) for the source and target of an arrow \(\edge\), and by convention we have \((\edge\rev)\rev=\edge\).

\begin{dfn}[\cite{gelfand_model_1979}]\label{dfn:ppa}
  Let \(\graph\) be a graph without loops and choose any map \(\map[\Qsgn]{\Edge{\double{\graph}}}{\set{-1, 1}}\) such that \(\Qsgn(\edge)\neq\Qsgn(\edge\rev)\) for all \(\edge\in\Edge{\double{\graph}}\). The \emph{preprojective algebra} of \((\graph, \Qsgn)\) is the unital associative algebra
  \begin{equation*}
    \Ppa \defn \quo*{\Comp\double{\graph}}{\paren*{\Sum{\Qsgn(\edge)\edge\edge\rev}[\edge\in\Edge{\double{\graph}}]}}.
  \end{equation*}
\end{dfn}

It has been shown (\cite[Lemma 2.2]{crawley-boevey_noncommutative_1998}) that this definition does not depend on the choice of map \(\Qsgn\), up to isomorphism. At each vertex \(\ver\in\Ver{\graph}\), there is a corresponding idempotent \(\idemp[\ver]\in\Ppa\), and multiplying both sides of the defining equation by each \(\idemp[\ver]\) gives the ideal
\begin{equation}\label{eqn:ppr}
  \idl*{\Ppr[\ver] \defn\Sum{\Qsgn(\edge)\edge\edge\rev}[\source[\edge]=\ver]}[\ver\in\Ver{\graph}],
\end{equation}
which equivalently presents the preprojective algebra.

\begin{thm}[{\cite[\S 6]{baer_preprojective_1987}}]\label{thm:ppa}
  If \(\graph\) is a connected graph without loops, then the following statements hold.
  \begin{enumerate}[label=(\arabic*)]
      \item \(\Ppa\) is finite-dimensional (thus noetherian) if and only if \(\graph\) is finite type.
      \item \(\Ppa\) is infinite-dimensional and noetherian if and only if \(\graph\) is affine type.
      \item \(\Ppa\) is not noetherian if and only if \(\graph\) is wild type.
  \end{enumerate}
\end{thm}

It follows that, in general, the category \(\fgMod{\Ppa}\) of finitely generated modules is not abelian, and so we must work with the full module category \(\Mod{\Ppa}\) in order to construct the derived category \(\Db{\Ppa}\). However, to remain close to the noetherian situation we make use of the fact that \(\Ppa\) admits an \(\Natural\)-grading by path length and consider \emph{nilpotent} modules. Denote by \(\rad\) the \emph{arrow ideal} in \(\Ppa\), the ideal generated by all length one paths in \(\Comp\double{\graph}\).

\begin{dfn}\label{dfn:nilp}
  A finite-dimensional (right) \(\Ppa\)-module \(\mod\) is said to be \emph{nilpotent} if there is some \(\len\in\Int[0]\) such that \(\mod\rad<\len>=0\).
\end{dfn}

We write \(\Nilp[\Ppa]\) for the category of nilpotent \(\Ppa\)-modules. Recall that the \emph{vertex simple} \(\simp[\ver]\) is the one-dimensional \(\Ppa\)-module that \(\idemp[\ver]\) acts on as the identity and all other arrows act as zero. It can be shown (see \cite[\S 2.4]{ringel_introduction_2012}) that \(\mod\in\Mod{\Ppa}\) is nilpotent precisely when it is finite-dimensional and all of its composition factors are vertex simples.

This means that every simple in \(\Nilp[\Ppa]\) is of the form \(\simp[\ver]\) for some \(\ver\in\Ver{\graph}\). By \cite[\S 4.2]{brenner_periodic_2002}, each vertex simple \(\simp[\ver]\) has an explicit projective resolution
\begin{equation}\label{eqn:projres}
  \begin{tikzcd}
    0 \ar[r] 
    & \proj[\ver] \ar[r]
    & \proj*[\ver] \defn \bigoplus\limits_{\substack{\edge*\in \Edge{\double{\graph}} \\\source[\edge*]=\ver}}\proj[\target[\edge*]] \ar[r] 
    & \proj[\ver] \ar[r, "{\mor[\ver]}"]
    & \simp[\ver] \ar[r]
    & 0,
  \end{tikzcd}
\end{equation}
where \(\proj[\ver] \defn \Ppa\idemp[\ver]\) is the \emph{vertex projective} and \(\map[\mor[\ver]]{\proj[\ver]}{\simp[\ver]}\) is the unique map in (graded) degree zero.

The category \(\Nilp[\Ppa]\) is a \emph{Serre subcategory} of \(\Mod{\Ppa}\), meaning it is closed under inclusions, quotients, and extensions. This allows us to define the following thick triangulated subcategory.

\begin{dfn}\label{dfn:gcat}
  For a connected graph \(\graph\) define the full subcategory
  \begin{equation*}
    \Gcat*[\graph]\defn\Db[\Nilp]{\Ppa} = \set{\obj\in\Db{\Ppa}}[\coH{\sdot}{\obj}\in\Nilp[\Ppa]].
  \end{equation*}
  of \(\Db{\Ppa}\) with nilpotent cohomology.
\end{dfn}

By \cite[\S 4]{keller_calabi-yau_2008}, \(\Gcat*[\graph]\) is a \(2\)-Calabi--Yau triangulated category (often termed a \emph{K3 category}) whenever \(\graph\) is not of finite type. This means that there is a functorial isomorphism
\begin{equation}\label{eqn:cy2}
  \map*{\Hom[\Gcat*[\graph]]{\obj}{\obj*}}{\Hom[\Gcat*[\graph]]{\obj*}{\obj\shift{2}}\dual}
\end{equation}
for all \(\obj,\obj*\in\Gcat*[\graph]\). As a consequence of \(\Ppa\) having global dimension two, we have
\begin{equation}\label{eqn:extfin}
  \Dim*[\Comp]{\direct*[\deg\in\Int]\Hom[\Gcat*[\graph]]{\obj}{\obj*\shift{\deg}}}< \infty,
\end{equation}
(see \cite[\S 1.1]{bridgeland_stability_2009}), a property often called \emph{finite type} in the literature, but which we call \emph{\(\mathrm{Ext}\)-finite} to avoid confusion with \say{finite type graph} used previously.

Now consider the K-theory of \(\Gcat*[\graph]\). First note that the category \(\Nilp[\Ppa]\) is a length category by definition, with finitely many isoclasses of simple objects \(\set{\simp[\ver]}[\ver\in\Ver{\graph}]\). Therefore, the Grothendieck group of \(\Gcat*[\graph]\) is spanned by the vertex simples. That is,
\begin{equation*}
  \Kgrp{\Gcat*[\graph]}=\Kgrp{\Nilp[\Ppa]}\iso\direct*[\ver\in\Ver{\graph}]\Int\Kthy{\simp[\ver]}.
\end{equation*}
The Ext-finite condition \eqref{eqn:extfin} allows us to equip \(\Kgrp{\Gcat*[\graph]}\) with the bilinear \emph{Euler form}
\begin{equation*}
  \map[\upchi]{\Kgrp{\Gcat*[\graph]}\times\Kgrp{\Gcat*[\graph]}}{\Int}, \qquad \paren{\Kthy{\obj}, \Kthy{\obj*}} \mapsto \Sum{(-1)^k\Dim[\Comp]{\Hom[\Gcat*[\graph]]{\obj}{\obj*\shift{\deg}}}}[\deg\in\Int].
\end{equation*}
This gives an identification with the root lattice, via
\begin{equation}\label{eqn:juniper}
  (\Kgrp{\Gcat*[\graph]}, \upchi) \cong (\Rts, \form{\blank}{\blank}),\qquad \Kthy{\simp[\ver]}\mapsto\rt[\ver].
\end{equation}

\begin{lem}\label{lem:nilpres}
  If \(\paren{\Db{\Ppa}\aisle, \Db{\Ppa}\aisle*}\) is the standard t-structure in \(\Db{\Ppa}\), which has heart \(\Mod{\Ppa}\), then
  \begin{equation*}
    \Gcat*[\graph]^{\,\leq0} \defn \Gcat*[\graph] \inter \Db{\Ppa}\aisle, \qquad \Gcat*[\graph]^{\,\geq0} \defn \Gcat*[\graph] \inter \Db{\Ppa}\aisle*
  \end{equation*}
  defines a t-structure in \(\Gcat*[\graph]\) with heart \(\Nilp[\Ppa]\). In other words, the standard t-structure in \(\Db{\Ppa}\) \emph{restricts} to the category \(\Gcat*[\graph]\).
\end{lem}
\begin{proof}
  Since \(\Nilp[\Ppa]\subseteq\Mod{\Ppa}\) is a Serre subcategory, this follows from \cite[\S 1.3.19]{beilinson_faisceaux_1982}.
\end{proof}

We call \(\Nilp[\Ppa]\) the \emph{standard heart} in \(\Gcat*[\graph]\). Using the tilting theory of \(\Mod{\Ppa}\) we can assign hearts to other alcoves in the Coxeter arrangement. To this end, let \(\ver\in\Ver{\graph}\) and consider the two-sided ideal
\begin{equation*}
  \tilt[\ver]\defn\Ppa(1-\idemp[\ver])\Ppa
\end{equation*}
in \(\Ppa\). Equivalently, this is the kernel of the projection \(\surj{\Ppa}{\simp[\ver]}\), which is spanned by every path except \(\idemp[\ver]\). By \cite[\S III.1]{buan_cluster_2009}, \(\tilt[\ver]\) is a classical tilting \(\Ppa\)-module with \(\End[\Ppa]{\tilt[\ver]}\iso\Ppa\). This means that the projective dimension of \(\tilt[\ver]\) is at most one, \(\Ext[\Ppa]{1}{\tilt[\ver]}{\tilt[\ver]} = 0\), and there exists a two-term resolution of \(\Ppa\) by modules in \(\add[\tilt[\ver]]\). Then \cite[\S 6]{rickard_morita_1989} implies that \(\Db{\Ppa}\) admits the derived autoequivalences
\begin{equation*}
  \smut[\ver]\defn\RHom{\Ppa}{\tilt[\ver]}{\blank} \qquad\qquad \smut[\ver]\inv = \LTen{\Ppa}{\blank}{\tilt[\ver]},
\end{equation*}
which we call \emph{mutation functors}. Furthermore, the product of \(\tilt[\ver]\) and \(\tilt[\ver*]\) is also classical tilting, so we can compose mutation functors using reduced expressions in the Weyl group of \(\graph\) (see \cite{buan_cluster_2009}). 

Now consider the \emph{dimension vector}
\begin{equation*}
  \map[\Kthy{\blank}]{\Db[\fdMod]{\Ppa}}{\Rts}, \qquad \mod\cx\mapsto\Sum{\paren{-1}<\deg>\Qdim{\coH{\deg}{\mod\cx}}}[\deg\in\Int],
\end{equation*}
named as such because \(\Kthy{\mod} = \Qdim{\mod}\) for all \(\mod\in\fdMod[\Ppa]\). Since \(\Nilp[\Ppa]\subset\fdMod[\Ppa]\), the function \(\Kthy{\blank}\) restricts to \(\Gcat*[\graph]\). Then \cite[Lemma 2.22 and Theorem 2.28]{sekiya_tilting_2013} imply that the diagram
\begin{equation}\label{eqn:wpair}
  \begin{tikzcd}[column sep = 5em]
    \Gcat*[\graph] \ar[r, leftrightarrow, "{\smut[\ver]}", "{\smut[\ver]\inv}"'] \ar[d, "\Kthy{\blank}"'] & 
        \Gcat*[\graph] \ar[d, "\Kthy{\blank}"] \\ 
        V \ar[r, leftrightarrow, "{\sref[\ver]}"'] & 
        V
    \end{tikzcd}
\end{equation}
commutes.

\begin{lem}\label{lem:smuts}
  If \(\ver,\ver*\in\Ver{\graph}\) are distinct, then the images of \(\simp[\ver]\) (viewed as a complex in degree zero) under the corresponding mutation functors satisfy
  \begin{equation*}
    \smut[\ver]<\pm>\simp[\ver]\iso\simp[\ver]\shift{\mp1}, \qquad \smut[\ver*]\simp[\ver]\iso\Hom[\Ppa]{\tilt[\ver*]}{\simp[\ver]}, \qquad \smut[\ver*]\inv\simp[\ver]=\simp[\ver]\tensor[\Ppa]\tilt[\ver*].
  \end{equation*}
  In particular, the last two modules have dimension vector \(\sref[\ver*]\rt[\ver]\).
\end{lem}
\begin{proof}
  This follows from the proof of \cite[Theorem 2.28]{sekiya_tilting_2013}, or equivalently \cite[Lemma 5.7]{wemyss_flops_2018}.
\end{proof}

\cref{lem:smuts} implies that the cohomology of the mutation functors applied to \(\nilp\in\Nilp[\Ppa]\) is nilpotent. Thus, \(\smut[\ver]<\pm>\) descends to \(\Gcat*[\graph]\) and generate a group of autoequivalences
\begin{equation*}
  \Br\Gcat*[\graph]\defn\gen{\smut[\ver]}[\ver\in\Ver{\graph}]\leq\Auteq{\Gcat*[\graph]}.
\end{equation*}
See \cite[\S 2.3]{sekiya_tilting_2013} for more details. By a proof analogous to that of \cite[Lemma 4.2]{donovan_stringy_2022}, if \(\func\in\Br\Gcat*[\graph]\) and \(\Abel*\defn\func(\Abel)\), then for all \(\obj\in\Gcat*[\graph]\) we have
\begin{equation}\label{eqn:mutcoh}
  \coH[\Abel*]{\sdot}{\func\obj} \iso \func\coH[\Abel]{\sdot}{\obj}.
\end{equation}

It is useful to also think of \(\smut[\ver]<\pm1>(\Nilp[\Ppa])\) in terms of the left/right simple tilts introduced in \cite{happel_tilting_1996}. Recall that a torsion pair in an abelian category \(\Abel\) is a pair of full subcategories \((\Tors, \Tors*)\), such that \(\Hom[\Abel]{\Tors}{\Tors*} = 0\) and every \(\mod\in\Abel\) fits into a short exact sequence of the form
\begin{equation}\label{eqn:tors}
  \begin{tikzcd}
    0 \ar[r] & \mod* \ar[r] & \mod \ar[r] & \mod*' \ar[r] & 0,
  \end{tikzcd}
\end{equation}
where \(\mod*\in\Tors\) and \(\mod*'\in\Tors*\). Now, each vertex simple \(\simp[\ver]\) induces the torsion theories \((\add[\simp[\ver]], \Tors*[\ver])\) and \((\Tors[\ver], \add[\simp[\ver]])\) in the standard heart \(\Nilp[\Ppa]\), where
\begin{equation*}
  \Tors[\ver] \defn \set{\nilp\in\Nilp[\Ppa]}[\Hom[\Nilp]{\nilp}{\simp[\ver]} = 0], \qquad \Tors*[\ver] \defn \set{\nilp\in\Nilp[\Ppa]}[\Hom[\Nilp]{\simp[\ver]}{\nilp} = 0].
\end{equation*}
Tilting \(\Nilp[\Ppa]\) with respect to these torsion theories yields, by \cite{happel_tilting_1996}, bounded hearts
\begin{align*}
  \Tilt[\ver]{\Nilp[\Ppa]} &\defn \set{\obj\in\Gcat*[\graph]}[\coH{0}{\obj}\in\Tors*[\ver], \coH{1}{\obj}\in\add[\simp[\ver]], \coH{\deg}{\obj} = 0 \text{ otherwise}], \\ 
  \Tilt*[\ver]{\Nilp[\Ppa]} &\defn \set{\obj\in\Gcat*[\graph]}[\coH{-1}{\obj}\in\add[\simp[\ver]], \coH{0}{\obj}\in\Tors[\ver], \coH{\deg}{\obj} = 0 \text{ otherwise}],
\end{align*}
which coincide with the image of the standard heart under mutation as follows.

\begin{lem}\label{lem:tiltmut}
  For all \(\ver\in\Ver{\graph}\) we have equalities
  \begin{equation*}
    \Tilt[\ver]{\Nilp[\Ppa]} = \smut[\ver](\Nilp[\Ppa]), \qquad \Tilt*[\ver]{\Nilp[\Ppa]} = \smut[\ver]\inv(\Nilp[\Ppa])
  \end{equation*}
  of bounded hearts in \(\Gcat*[\graph]\).
\end{lem}
\begin{proof}
  This is similar to \cite[\S 4.2]{donovan_stringy_2022}, but we give the idea in our context. Since every abelian category in the statement is a bounded heart, it suffices to show that \(\smut[\ver]\inv(\Nilp[\Ppa]) \subseteq \Tilt*[\ver]{\Nilp[\Ppa]}\) and similarly for the left tilt. This is because of the well-known (see, for example, \cite[Exercise 5.6]{macri_lectures_2017}) fact that if \(\Abel*\subseteq\Abel\) are bounded hearts in a triangulated category \(\Gcat\), then \(\Abel* = \Abel\). As every \(\nilp\in\Nilp[\Ppa]\) is built out of vertex simples, we just check \(\smut[\ver]\inv(\simp[\ver]),\smut[\ver]\inv(\simp[\ver*])\in\Tilt*[\ver]{\Nilp[\Ppa]}\). By \cref{lem:smuts}, \(\smut[\ver]\inv(\simp[\ver])\) and \(\smut[\ver]\inv(\simp[\ver*])\) are in cohomological degree \(-1\) and \(0\) respectively. Firstly
  \begin{equation*}
    \coH{-1}{\smut[\ver]\inv(\simp[\ver])} = \coH{-1}{\simp[\ver]\shift{1}} = \simp[\ver] \in \add[\simp[\ver]],
  \end{equation*}
  as required. Then adjunction, along with \cref{lem:smuts}, implies
  \begin{equation*}
    \Hom[\Nilp]{\coH{0}{\smut[\ver]\inv\simp[\ver*]}}{\simp[\ver]} \iso \Hom[\Gcat]{\smut[\ver]\inv\simp[\ver*]}{\simp[\ver]} \iso \Hom[\Gcat]{\simp[\ver*]}{\smut[\ver]\simp[\ver]} = \Hom[\Gcat]{\simp[\ver*]}{\simp[\ver]\shift{-1}} = 0.
  \end{equation*}
  That is, \(\coH{0}{\smut[\ver]\inv\simp[\ver*]}\) is torsion. This implies that \(\smut[\ver]\inv(\simp[\ver*])\in\Tilt*[\ver]{\Nilp[\Ppa]}\).
\end{proof}

This gives us a way of getting new hearts in \(\Gcat*[\graph]\) by applying the mutation functors repeatedly to \(\Nilp[\Ppa]\).

  \section{Real charges and real flows}\label{sec:real}
  In this section we equip the pair \((\Hpln[\graph],\Gcat*[\graph])\) constructed in \cref{sec:back} with a real central charge and an assignment of hearts to alcoves. Roughly speaking, the former consists of the embedding of \(\Level[\graph]\) inside \(\Ths\), and the latter takes the alcove determined by \(\weyl\in\Weyl\) and assigns to it the image of \(\Nilp[\Ppa]\) under some composition \(\smut[\ver]<\pm>\smut[\ver*]<\pm>\cdots\) of mutation functors. The choice of signs in this assignment naturally leads to the definition of a real flow.

\subsection{The central charge}\label{ssec:charge}

Let \(\inj[\charge]{\Level[\graph]}{\Ths}\) be the inclusion map. Then, by \eqref{eqn:juniper}, \(\charge\) is a map to \(\Kgrp{\Gcat*[\graph]}[\Real]\dual\), as required by \cref{dfn:realstab}. Since \(\map[\charge]{\Level}{\Ths}\) is injective and \(\sref[\ver]\image{\charge} \subseteq \image{\charge}\), there is an action of \(\Weyl[\graph]\) on \(\Level\) making the diagram
\begin{equation}\label{eqn:charge}
  \begin{tikzcd}
    \Level[\graph]\ar[r, "\charge", hook] \ar[d, "{\sref[\ver]}"'] & \Ths \ar[r, "\Kthy{\mod}"] \ar[d, "{\sref[\ver]}"] & \Real \ar[d, equal] \\
    \Level[\graph]\ar[r, "\charge"', hook] & \Ths \ar[r, "{[\smut[\ver]\mod]}"'] & \Real
  \end{tikzcd}
\end{equation}
commute for all \(\mod\in\Gcat*[\graph]\). By parametrising \(\Level[\graph]\) as (\(\rank-1\))-dimensional affine or hyperbolic space, we can give an explicit formula for the real central charge \(\charge\). This is a real polynomial for \(\graph\) affine, in line with the original definition of \cite{anno_stability_2015}, but becomes a real power series when \(\graph\) is hyperbolic.

\begin{lem}\label{lem:param}
  For \(\graph\) hyperbolic and \(\stab*\in\Level[\graph]\), \(\evalue*[\ver]\in\Real\minus\set{0}\) and \(1\leq \ver \leq \rank\) define
  \begin{equation*}\label{eqn:param}
    \stab[\ver] \defn \frac{1}{\evalue*[\ver]}\sinh{\stab*[\rank-\ver]}\Prod{\cosh{\stab*[\rank-\ver*]}}[\ver*=1]<\ver-1>,
  \end{equation*}
  with conventions that the empty product and \(\sinh{\stab*[0]}\) are \(1\). Then
  \begin{equation*}
    -\evalue*[\rank]<2>\stab[\rank]<2> + \Sum{\evalue*[\ver]<2>\stab[\ver]<2>}[\ver=1]<\rank-1> = -1,
  \end{equation*}
  and hence we have a parametrisation of the hyperboloid in \(\Real<\rank>\). 
\end{lem}
\begin{proof}
  This is clear by using the identity \(\cosh^2{\stab*[\ver]} - \sinh^2{\stab*[\ver]} = 1\).
\end{proof}

If we take \(\evalue*[\ver] = \sqrt{\evalue[\ver]}\) to be the eigenvalues of \(\GCM[\graph]\), then the above is a parametrisation of \(\Level\).

\begin{exm}\label{exm:param}
  Consider \(\graph = \Kronecker[\val]\) from \cref{exm:rktwo}, with \(\val\geq3\). As \(\Level[\graph]\iso\Real\), we may parametrise the real central charge as follows
  \begin{equation*}
    \map[\charge[\val]]{\Real}{\Ths}, \qquad\real\mapsto\frac{2}{\sqrt{\val<2>-4}}\Matrix{\sinh\log\Kcon[\val]<1-\real>\\\sinh\log\Kcon[\val]<\real>}=\frac{1}{\Kcon[\val]-\Kcon[\val]\inv}\Matrix{\Kcon[\val]<1-\real>-\Kcon[\val]<\real-1>\\\Kcon[\val]<\real>-\Kcon[\val]<-\real>},
  \end{equation*}
  recalling the Kronecker constant \(\Kcon[\val]\) from \cref{exm:rktwo}. This satisfies 
  \begin{equation*}
    \charge[\val]\paren{\real}[\xver]<2> + \val\charge[\val]\paren{\real}[\xver]\charge[\val]\paren{\real}[1] + \charge[\val]\paren{\real}[1]<2> = 1
  \end{equation*}
  for all \(\real\in\Real\) and can be compared with Figure 1. Viewing \(\val\) as being continuous, taking the limit as \(\val\) approaches \(2\) recovers the classical affine level. Indeed, a simple application of \Hopital* rule yields
  \begin{equation*}
    \charge[2](\real) = \Lim{\charge[\val](\real)}[\val\to 2] = \paren{1-\real, \real},
  \end{equation*}
  which is Figure 1(2).
\end{exm}

\subsection{Real functors and real flows}\label{ssec:flow}

Below, the assignment of t-structures to alcoves is controlled by a functor \(\flow\) from the arrangement groupoid to the autoequivalence groupoid, called a \emph{real flow}. Essentially, this is a choice of \(\above\) or \(\below\) between neighbouring alcoves, categorifying the idea of a \say{direction} on \(\Alc\). Recall from \cref{ssec:cox} the free category (or Deligne quiver) \(\Free\Hpln[\graph]\) of the Coxeter arrangement associated with \(\graph\). We may think of each wall-crossing as either a length one reduced positive path or a simple reflection via the identification \(\alc=\weyl\fAlc, \alc*=\weyl\sref[\ver]\weyl\inv\alc\). The common abuse of notion \(\sref[\ver]\) will be used for both directions of the wall-crossing, as this becomes important when we study equalities in the braid group.

\begin{dfn}\label{dfn:func}
  A \emph{real functor} on \(\Hpln[\graph]\) is a functor 
  \begin{equation*}
    \map[\flow]{\Free\Hpln[\graph]}{\gAuteq{\Gcat*[\graph]}}, \qquad \path\mapsto\flow[\path],
  \end{equation*}
  sending every object \(\alc\) to \(\Gcat*[\graph]\) and every generating morphism \(\map[\sref[\ver]]{\alc}{\alc*}\) to either \(\smut[\ver]\) or \(\smut[\ver]\inv\), such that \(\flow[\map{\alc*}{\alc}] = \paren{\flow[\map{\alc}{\alc*}]}\inv\). Here \(\ver\in\Ver{\graph}\) corresponds to the wall between \(\alc\) and \(\alc*\), in the sense that if \(\stab\in\alc=\weyl\fAlc\) for some \(\weyl\in\Weyl\), then \(\weyl\sref[\ver]\weyl\inv\stab\in\alc*\). If \(\flow[\map{\alc*}{\alc}] = \smut[\ver]\), then we say that \(\alc\) is \emph{above} \(\alc*\) and write \(\alc\above\alc*\), otherwise \(\alc\) is \emph{below} \(\alc*\), written \(\alc\below\alc*\). We can summarise this as
  \begin{equation*}
    \begin{tikzcd}
      \alc \ar[r, shift left, "{\sref[\ver]}"] & \alc* \ar[l, shift left, "{\sref[\ver]}"]
    \end{tikzcd} \qquad \xmapsto{\phantom{m}\flow\phantom{m}} \qquad
    \begin{cases}
      \begin{tikzcd}
        \Gcat*[\graph] \ar[r, shift left, "{\smut[\ver]}"] & \Gcat*[\graph] \ar[l, shift left, "{\smut[\ver]\inv}"]
      \end{tikzcd}
       & \text{if } \alc\above\alc*, \\
      \begin{tikzcd}
        \Gcat*[\graph] \ar[r, shift left, "{\smut[\ver]\inv}"] & \Gcat*[\graph] \ar[l, shift left, "{\smut[\ver]}"]
      \end{tikzcd}
      & \text{if } \alc\below\alc*.
    \end{cases}
  \end{equation*}
  If, around every flat in \(\Hpln[\graph]\), there is a unique source alcove \(\alc[\source]\), opposite which is a unique target alcove \(\alc[\target]\), then \(\flow\) is a \emph{real flow}.
\end{dfn}

The use of the term \say{real} in this definition is justified in \cref{lem:flowtriv}. In the original definition, \cite{anno_stability_2015} work in the affine setting and define a direction on \(\Alc\) by first linearising every hyperplane in \(\Hpln\) and then choosing a connected component of the finite Coxeter arrangement to act as a \say{positive} direction. With no reasonable analogue in the hyperbolic setting, \cref{dfn:func} above categorifies this idea and, in turn, removes obstructions to finding real variations of stability in general.

\begin{rmk}\label{rmk:func}
  The choice of the \say{above} alcove as being the inverse mutation functor (or equivalently, by \cref{lem:tiltmut}, the right tilt) is in light of \cref{lem:smuts}, since we wish to line up the heart assignment with the shift in \eqref{eqn:quoprop}. Such an inverse arising is common, for example in the Bridgeland--Chen flop functor. However, it is important to ensure that the heart assignment induced by \cref{dfn:func} is well-defined, irrespective of minimal expressions for the Weyl element connecting two alcoves.
\end{rmk}

\begin{exm}\label{exm:bruhat}
  One example of a real flow comes from the \emph{Bruhat order} \(\Bruhat\) on \(\Weyl[\graph]\), where we define \(\weyl\fAlc = \alc \below \alc* = \weyl\sref[\ver]\fAlc\) if and only if \(\weyl\Bruhat\weyl\sref[\ver]\). To see this, fix a flat of type \(\finA[2]\) (the case \(\finA[1]\times\finA[1]\) is easier) and consider the three adjacent hyperplanes \(\hpln[1],\hpln[2],\hpln[3]\), viewed in \(\TC\). The alcove \(\fAlc\) is in one of the six regions in \(\TC\) that form the connected components of the complement of these these hyperplanes. Each of the six alcoves neighbouring the flat sits in precisely one of these regions, and therefore the alcove \(\alc\) sitting in the same region as \(\fAlc\) is the smallest with respect to \(\Bruhat\), since any reduced path \(\map{\fAlc}{\alc}\) does not cross any of \(\hpln[1],\hpln[2],\hpln[3]\).
  \begin{equation*}
    \begin{tikzpicture}[scale=2, decoration = snake]
      \node (Af) at (-2, -0.75) {\(\fAlc\)};
      \node (A) at (-{sqrt(3)/3}, -1/3)  {\(\alc\)};
      \node[label={right:\(\hpln[3]\)}] (H3) at (1/2, -{sqrt(3)/2}) {};
      \node[label={right:\(\hpln[2]\)}] (H2) at (1, 0) {};
      \node[label={right:\(\hpln[1]\)}] (H1) at (1/2, {sqrt(3)/2}) {};
      \node[circle, fill, scale=.5] (F) at (0, 0) {};
      \tkzDefPoint(0,0){O}
      \tkzDefPoint(1,0){P}
      \hgline{0}{180}
      \hgline{60}{240}
      \hgline{120}{300}
      \draw[color=black, decorate, ->] (Af) -- (A) node[above, midway]{\(\weyl\)};
    \end{tikzpicture}
  \end{equation*}
  We next claim that \(\alc\) is the unique target alcove in the flat. The claim holds since the path \(\weyl\sref[\ver]\) is reduced (crossing \(\hpln[1]\) only once), as is \(\weyl\sref[\ver]\sref[\ver*]\) (crossing \(\hpln[1]\) and \(\hpln[2]\) only once) and \(\weyl\sref[\ver]\sref[\ver*]\sref[\ver]\). Thus \(\weyl\Bruhat\weyl\sref[\ver]\Bruhat\weyl\sref[\ver]\sref[\ver*]\Bruhat\weyl\sref[\ver]\sref[\ver*]\sref[\ver]\). The other direction, namely \(\weyl\Bruhat\weyl\sref[\ver*]\Bruhat\weyl\sref[\ver*]\sref[\ver]\Bruhat\weyl\sref[\ver*]\sref[\ver]\sref[\ver*]\), is similar.
\end{exm}

Since the word problem is solved for the braid group \(\Br[\finA[2]]\), we may easily verify the following cases.

\begin{lem}\label{lem:flow}
  Let \(\graph\) be \(\finA[1]\times\finA[1]\) or \(\finA[2]\), and let \(\ver,\ver*\in\Ver{\graph}\) be distinct. The following are the only braid equalities of the form \(\sref[\ver]<\pm>\sref[\ver*]<\pm> = \sref[\ver*]<\pm>\sref[\ver]<\pm>\), respectively \(\sref[\ver]<\pm>\sref[\ver*]<\pm>\sref[\ver]<\pm> = \sref[\ver*]<\pm>\sref[\ver]<\pm>\sref[\ver*]<\pm>\), that hold in \(\Br[\graph]\).
  \begin{align*}
    \finA[1]&\times\finA[1] & &\finA[2] \\
    \sref[\ver]\sref[\ver*] &= \sref[\ver*]\sref[\ver] & \sref[\ver]\sref[\ver*]\sref[\ver] &= \sref[\ver*]\sref[\ver]\sref[\ver*]  \\
    \sref[\ver]\sref[\ver*]\inv &= \sref[\ver*]\inv\sref[\ver] & \sref[\ver]\sref[\ver*]\sref[\ver]\inv &= \sref[\ver*]\inv\sref[\ver]\sref[\ver*] \\
    \sref[\ver]\inv\sref[\ver*] &= \sref[\ver*]\sref[\ver]\inv & \sref[\ver]\inv\sref[\ver*]\sref[\ver] &= \sref[\ver*]\sref[\ver]\sref[\ver*]\inv \\
    \sref[\ver]\inv\sref[\ver*]\inv &= \sref[\ver*]\inv\sref[\ver]\inv & \sref[\ver]\sref[\ver*]\inv\sref[\ver]\inv &= \sref[\ver*]\inv\sref[\ver]\inv\sref[\ver*] \\
    & & \sref[\ver]\inv\sref[\ver*]\inv\sref[\ver] &= \sref[\ver*]\sref[\ver]\inv\sref[\ver*]\inv \\
    & & \sref[\ver]\inv\sref[\ver*]\inv\sref[\ver]\inv &= \sref[\ver*]\inv\sref[\ver]\inv\sref[\ver*]\inv
  \end{align*}
  Each corresponds to a longest diagonal across a \(4\)-gon and a \(6\)-gon, respectively.
\end{lem}
\begin{proof}
  This is readily checked, for example, by the computer algebra package \texttt{Magma} \cite{bosma_magma_1997}, iterating over the \(2^6 = 64\) choices of positive and negative powers in relations of the form \(\sref[\ver]<\pm>\sref[\ver*]<\pm>\sref[\ver]<\pm> = \sref[\ver*]<\pm>\sref[\ver]<\pm>\sref[\ver*]<\pm>\). Each of these expressions determines a source and target alcove. Indeed, fixing a reference alcove \(\alc\), labelling walls with alternating indices, and thinking of wall-crossing across \(\ver\) to a higher alcove as multiplication by \(\sref[\ver]\), we obtain the stated expressions by passing from \(\alc\) to its opposite \(\cl{\alc}\) in both ways. The \(\finA[1]\times\finA[1]\) case is similar.
  \begin{equation}\label{eqn:flow}
    \begin{array}{ccc}
      \begin{tikzpicture}[scale=1.5]
        \node[right] (z1) at (1, 0) {\(i\)};
        \node[above] (z2) at (1/2, {sqrt(3)/2}) {\(j\)};
        \node[above] (z3) at (-1/2, {sqrt(3)/2}) {\(i\)};
        \node[left] (z4) at (-1, 0) {\(j\)};
        \node[below] (z5) at (-1/2, -{sqrt(3)/2}) {\(i\)};
        \node[below] (z6) at (1/2, -{sqrt(3)/2}) {\(j\)};
        \node (A) at (0, -1/3) {\(\alc\)};
        \node (Aop) at (0, 1/3) {\(\cl{\alc}\)};
        \node[r, rotate=90] (p1) at ({2*sqrt(3)/9}, 0) {\Large\(\prec\)};
        \node[r, rotate=150] (p2) at ({sqrt(3)/9}, 1/3) {\Large\(\prec\)};
        \node[r, rotate=30] (p3) at (-{sqrt(3)/9}, 1/3) {\Large\(\prec\)};
        \node[r, rotate=90] (p4) at (-{2*sqrt(3)/9}, 0) {\Large\(\prec\)};
        \node[r, rotate=150] (p5) at (-{sqrt(3)/9}, -1/3) {\Large\(\prec\)};
        \node[r, rotate=30] (p6) at ({sqrt(3)/9}, -1/3) {\Large\(\prec\)};
        \node (A12) at ({sqrt(3)/8}, 1/8) {};
        \node (A23) at (0, 1/4) {};
        \node (A34) at (-{sqrt(3)/8}, 1/8) {};
        \node (A45) at (-{sqrt(3)/8}, -1/8) {};
        \node (A56) at (0, -1/4) {};
        \node (A61) at ({sqrt(3)/8}, -1/8) {};
        \tkzDefPoint(0,0){O}
        \tkzDefPoint(1,0){P}
        \tkzClipCircle(O,P)
        \hgline{0}{60}
        \hgline{60}{120}
        \hgline{120}{180}
        \hgline{180}{240}
        \hgline{240}{300}
        \hgline{300}{360}
        \hgline{0}{180}
        \hgline{60}{240}
        \hgline{120}{300}
        \draw[->, r, thick] (A23) edge (A56);
      \end{tikzpicture}
      &
      \begin{tikzpicture}[scale=1.5]
        \node[right] (z1) at (1, 0) {\(i\)};
        \node[above] (z2) at (1/2, {sqrt(3)/2}) {\(j\)};
        \node[above] (z3) at (-1/2, {sqrt(3)/2}) {\(i\)};
        \node[left] (z4) at (-1, 0) {\(j\)};
        \node[below] (z5) at (-1/2, -{sqrt(3)/2}) {\(i\)};
        \node[below] (z6) at (1/2, -{sqrt(3)/2}) {\(j\)};
        \node (A) at (0, -1/3) {\(\alc\)};
        \node (Aop) at (0, 1/3) {\(\cl{\alc}\)};
        \node[r, rotate=90] (p1) at ({2*sqrt(3)/9}, 0) {\Large\(\prec\)};
        \node[r, rotate=150] (p2) at ({sqrt(3)/9}, 1/3) {\Large\(\prec\)};
        \node[r, rotate=210] (p3) at (-{sqrt(3)/9}, 1/3) {\Large\(\prec\)};
        \node[r, rotate=90] (p4) at (-{2*sqrt(3)/9}, 0) {\Large\(\prec\)};
        \node[r, rotate=150] (p5) at (-{sqrt(3)/9}, -1/3) {\Large\(\prec\)};
        \node[r, rotate=210] (p6) at ({sqrt(3)/9}, -1/3) {\Large\(\prec\)};
        \node (A12) at ({sqrt(3)/8}, 1/8) {};
        \node (A23) at (0, 1/4) {};
        \node (A34) at (-{sqrt(3)/8}, 1/8) {};
        \node (A45) at (-{sqrt(3)/8}, -1/8) {};
        \node (A56) at (0, -1/4) {};
        \node (A61) at ({sqrt(3)/8}, -1/8) {};
        \tkzDefPoint(0,0){O}
        \tkzDefPoint(1,0){P}
        \tkzClipCircle(O,P)
        \hgline{0}{60}
        \hgline{60}{120}
        \hgline{120}{180}
        \hgline{180}{240}
        \hgline{240}{300}
        \hgline{300}{360}
        \hgline{0}{180}
        \hgline{60}{240}
        \hgline{120}{300}
        \draw[->, r, thick] (A34) edge (A61);
      \end{tikzpicture}
      &
      \begin{tikzpicture}[scale=1.5]
        \node[right] (z1) at (1, 0) {\(i\)};
        \node[above] (z2) at (1/2, {sqrt(3)/2}) {\(j\)};
        \node[above] (z3) at (-1/2, {sqrt(3)/2}) {\(i\)};
        \node[left] (z4) at (-1, 0) {\(j\)};
        \node[below] (z5) at (-1/2, -{sqrt(3)/2}) {\(i\)};
        \node[below] (z6) at (1/2, -{sqrt(3)/2}) {\(j\)};
        \node (A) at (0, -1/3) {\(\alc\)};
        \node (Aop) at (0, 1/3) {\(\cl{\alc}\)};
        \node[r, rotate=270] (p1) at ({2*sqrt(3)/9}, 0) {\Large\(\prec\)};
        \node[r, rotate=150] (p2) at ({sqrt(3)/9}, 1/3) {\Large\(\prec\)};
        \node[r, rotate=210] (p3) at (-{sqrt(3)/9}, 1/3) {\Large\(\prec\)};
        \node[r, rotate=270] (p4) at (-{2*sqrt(3)/9}, 0) {\Large\(\prec\)};
        \node[r, rotate=150] (p5) at (-{sqrt(3)/9}, -1/3) {\Large\(\prec\)};
        \node[r, rotate=210] (p6) at ({sqrt(3)/9}, -1/3) {\Large\(\prec\)};
        \node (A12) at ({sqrt(3)/8}, 1/8) {};
        \node (A23) at (0, 1/4) {};
        \node (A34) at (-{sqrt(3)/8}, 1/8) {};
        \node (A45) at (-{sqrt(3)/8}, -1/8) {};
        \node (A56) at (0, -1/4) {};
        \node (A61) at ({sqrt(3)/8}, -1/8) {};
        \tkzDefPoint(0,0){O}
        \tkzDefPoint(1,0){P}
        \tkzClipCircle(O,P)
        \hgline{0}{60}
        \hgline{60}{120}
        \hgline{120}{180}
        \hgline{180}{240}
        \hgline{240}{300}
        \hgline{300}{360}
        \hgline{0}{180}
        \hgline{60}{240}
        \hgline{120}{300}
        \draw[->, r, thick] (A45) edge (A12);
      \end{tikzpicture}
      \\
      (---) = (---) & (--+) = (+--) & (-++) = (++-)
      \\\\
      \begin{tikzpicture}[scale=1.5]
        \node[right] (z1) at (1, 0) {\(i\)};
        \node[above] (z2) at (1/2, {sqrt(3)/2}) {\(j\)};
        \node[above] (z3) at (-1/2, {sqrt(3)/2}) {\(i\)};
        \node[left] (z4) at (-1, 0) {\(j\)};
        \node[below] (z5) at (-1/2, -{sqrt(3)/2}) {\(i\)};
        \node[below] (z6) at (1/2, -{sqrt(3)/2}) {\(j\)};
        \node (A) at (0, -1/3) {\(\alc\)};
        \node (Aop) at (0, 1/3) {\(\cl{\alc}\)};
        \node[r, rotate=270] (p1) at ({2*sqrt(3)/9}, 0) {\Large\(\prec\)};
        \node[r, rotate=330] (p2) at ({sqrt(3)/9}, 1/3) {\Large\(\prec\)};
        \node[r, rotate=210] (p3) at (-{sqrt(3)/9}, 1/3) {\Large\(\prec\)};
        \node[r, rotate=270] (p4) at (-{2*sqrt(3)/9}, 0) {\Large\(\prec\)};
        \node[r, rotate=330] (p5) at (-{sqrt(3)/9}, -1/3) {\Large\(\prec\)};
        \node[r, rotate=210] (p6) at ({sqrt(3)/9}, -1/3) {\Large\(\prec\)};
        \node (A12) at ({sqrt(3)/8}, 1/8) {};
        \node (A23) at (0, 1/4) {};
        \node (A34) at (-{sqrt(3)/8}, 1/8) {};
        \node (A45) at (-{sqrt(3)/8}, -1/8) {};
        \node (A56) at (0, -1/4) {};
        \node (A61) at ({sqrt(3)/8}, -1/8) {};
        \tkzDefPoint(0,0){O}
        \tkzDefPoint(1,0){P}
        \tkzClipCircle(O,P)
        \hgline{0}{60}
        \hgline{60}{120}
        \hgline{120}{180}
        \hgline{180}{240}
        \hgline{240}{300}
        \hgline{300}{360}
        \hgline{0}{180}
        \hgline{60}{240}
        \hgline{120}{300}
        \draw[->, r, thick] (A56) edge (A23);
      \end{tikzpicture}
      &
      \begin{tikzpicture}[scale=1.5]
        \node[right] (z1) at (1, 0) {\(i\)};
        \node[above] (z2) at (1/2, {sqrt(3)/2}) {\(j\)};
        \node[above] (z3) at (-1/2, {sqrt(3)/2}) {\(i\)};
        \node[left] (z4) at (-1, 0) {\(j\)};
        \node[below] (z5) at (-1/2, -{sqrt(3)/2}) {\(i\)};
        \node[below] (z6) at (1/2, -{sqrt(3)/2}) {\(j\)};
        \node (A) at (0, -1/3) {\(\alc\)};
        \node (Aop) at (0, 1/3) {\(\cl{\alc}\)};
        \node[r, rotate=270] (p1) at ({2*sqrt(3)/9}, 0) {\Large\(\prec\)};
        \node[r, rotate=330] (p2) at ({sqrt(3)/9}, 1/3) {\Large\(\prec\)};
        \node[r, rotate=30] (p3) at (-{sqrt(3)/9}, 1/3) {\Large\(\prec\)};
        \node[r, rotate=270] (p4) at (-{2*sqrt(3)/9}, 0) {\Large\(\prec\)};
        \node[r, rotate=330] (p5) at (-{sqrt(3)/9}, -1/3) {\Large\(\prec\)};
        \node[r, rotate=30] (p6) at ({sqrt(3)/9}, -1/3) {\Large\(\prec\)};
        \node (A12) at ({sqrt(3)/8}, 1/8) {};
        \node (A23) at (0, 1/4) {};
        \node (A34) at (-{sqrt(3)/8}, 1/8) {};
        \node (A45) at (-{sqrt(3)/8}, -1/8) {};
        \node (A56) at (0, -1/4) {};
        \node (A61) at ({sqrt(3)/8}, -1/8) {};
        \tkzDefPoint(0,0){O}
        \tkzDefPoint(1,0){P}
        \tkzClipCircle(O,P)
        \hgline{0}{60}
        \hgline{60}{120}
        \hgline{120}{180}
        \hgline{180}{240}
        \hgline{240}{300}
        \hgline{300}{360}
        \hgline{0}{180}
        \hgline{60}{240}
        \hgline{120}{300}
        \draw[->, r, thick] (A61) edge (A34);
      \end{tikzpicture}
      &
      \begin{tikzpicture}[scale=1.5]
        \node[right] (z1) at (1, 0) {\(i\)};
        \node[above] (z2) at (1/2, {sqrt(3)/2}) {\(j\)};
        \node[above] (z3) at (-1/2, {sqrt(3)/2}) {\(i\)};
        \node[left] (z4) at (-1, 0) {\(j\)};
        \node[below] (z5) at (-1/2, -{sqrt(3)/2}) {\(i\)};
        \node[below] (z6) at (1/2, -{sqrt(3)/2}) {\(j\)};
        \node (A) at (0, -1/3) {\(\alc\)};
        \node (Aop) at (0, 1/3) {\(\cl{\alc}\)};
        \node[r, rotate=90] (p1) at ({2*sqrt(3)/9}, 0) {\Large\(\prec\)};
        \node[r, rotate=330] (p2) at ({sqrt(3)/9}, 1/3) {\Large\(\prec\)};
        \node[r, rotate=30] (p3) at (-{sqrt(3)/9}, 1/3) {\Large\(\prec\)};
        \node[r, rotate=90] (p4) at (-{2*sqrt(3)/9}, 0) {\Large\(\prec\)};
        \node[r, rotate=330] (p5) at (-{sqrt(3)/9}, -1/3) {\Large\(\prec\)};
        \node[r, rotate=30] (p6) at ({sqrt(3)/9}, -1/3) {\Large\(\prec\)};
        \node (A12) at ({sqrt(3)/8}, 1/8) {};
        \node (A23) at (0, 1/4) {};
        \node (A34) at (-{sqrt(3)/8}, 1/8) {};
        \node (A45) at (-{sqrt(3)/8}, -1/8) {};
        \node (A56) at (0, -1/4) {};
        \node (A61) at ({sqrt(3)/8}, -1/8) {};
        \tkzDefPoint(0,0){O}
        \tkzDefPoint(1,0){P}
        \tkzClipCircle(O,P)
        \hgline{0}{60}
        \hgline{60}{120}
        \hgline{120}{180}
        \hgline{180}{240}
        \hgline{240}{300}
        \hgline{300}{360}
        \hgline{0}{180}
        \hgline{60}{240}
        \hgline{120}{300}
        \draw[->, r, thick] (A12) edge (A45);
      \end{tikzpicture}
      \\
      (+++) = (+++) & (++-) = (-++) & (++-) = (-++)
    \end{array}
  \end{equation}
\end{proof}

We are now ready to prove the main result of this section, which is that a real flow is precisely the information we need to get a well-defined assignment of hearts to alcoves.

\begin{ppn}\label{ppn:flow}
  If \(\flow\) is a real functor on \(\Hpln[\graph]\), then \(\flow\) is a flow if and only if it descends to a well-defined functor \(\map[\flow]{\gDel[\graph]}{\gAuteq{\Gcat*[\graph]}}\) on the arrangement groupoid. In this case, \(\flow[\path] = \flow[\path*]\) whenever \(\path,\path*\) are paths around a codimension two wall from an alcove to its opposite.
\end{ppn}
\begin{proof}
  Consider \(\Path\Hpln[\graph] = \quo{\Free\Hpln[\graph]}{\!\sim}\) defined in \cref{ssec:gdel}. If \(\sim_M\) is the smallest equivalence relation in \(\Free\Hpln[\graph]\) generated by codimension two (that is, braid) relations, then Matsumoto's \cref{thm:matsu} implies that \(\Path\Hpln[\graph] = \quo{\Free\Hpln[\graph]}{\!\sim_M}\). This is because any path in the free category can be identified with an element \(\weyl\) of the Weyl group, and any two reduced paths give rise to reduced expressions for \(\weyl\) and hence differ by a sequence of braid moves. Hence it suffices to check that the paths around codimension two walls give equivalent composition of autoequivalences. 
  
  First suppose that \(\flow\) is a real flow, so (possibly after rotating) each flat locally has the form
  \begin{equation*}
    \begin{array}{cc}
      \begin{tikzpicture}[scale=2]
        \node (A) at (0, -1/3) {\(\alc\)};
        \node (Aop) at (0, 1/3) {\(\cl{\alc}\)};
        \node[circle, fill, scale=.5] (F) at (0, 0) {};
        \tkzDefPoint(0,0){O}
        \tkzDefPoint(1,0){P}
        \tkzClipCircle(O,P)
        \hgline{0}{60}
        \hgline{60}{120}
        \hgline{120}{180}
        \hgline{180}{240}
        \hgline{240}{300}
        \hgline{300}{360}
        \hgline{0}{180}
        \hgline{60}{240}
        \hgline{120}{300}
        \node[r, rotate=90] (p1) at ({2*sqrt(3)/9}, 0) {\LARGE\(\prec\)};
        \node[r, rotate=150] (p2) at ({sqrt(3)/9}, 1/3) {\LARGE\(\prec\)};
        \node[r, rotate=30] (p3) at (-{sqrt(3)/9}, 1/3) {\LARGE\(\prec\)};
        \node[r, rotate=90] (p4) at (-{2*sqrt(3)/9}, 0) {\LARGE\(\prec\)};
        \node[r, rotate=150] (p5) at (-{sqrt(3)/9}, -1/3) {\LARGE\(\prec\)};
        \node[r, rotate=30] (p6) at ({sqrt(3)/9}, -1/3) {\LARGE\(\prec\)};
      \end{tikzpicture}
      &
      \begin{tikzpicture}[scale=2]
        \node (A) at (0, -1/4) {\(\alc\)};
        \node (Aop) at (0, 1/4) {\(\cl{\alc}\)};
        \node[circle, fill, scale=.5] (F) at (0, 0) {};
        \tkzDefPoint(0,0){O}
        \tkzDefPoint(1,0){P}
        \tkzClipCircle(O,P)
        \hgline{45}{135}
        \hgline{135}{225}
        \hgline{225}{315}
        \hgline{315}{45}
        \hgline{45}{225}
        \hgline{135}{315}
        \node[r, rotate=135] (p1) at (.25, .25) {\LARGE\(\prec\)};
        \node[r, rotate=45] (p2) at (-.25, .25) {\LARGE\(\prec\)};
        \node[r, rotate=135] (p3) at (-.25, -.25) {\LARGE\(\prec\)};
        \node[r, rotate=45] (p4) at (.25, -.25) {\LARGE\(\prec\)};
      \end{tikzpicture}
    \end{array}
  \end{equation*}
  Now by \cref{lem:flow} each of the six (or four) relations in \(\Path\Hpln[\graph]\) around this flat hold in \(\Br[\finA[2]]\) (or \(\Br[\finA[1]\times\finA[1]]\)), so by \cite[Proposition 2.27]{sekiya_tilting_2013} applied to \cite[III.1.8]{buan_cluster_2009}, the relations also hold in \(\gAuteq{\Gcat*[\graph]}\). Thus, \(\flow\) descends.

  Conversely, suppose that \(\flow\) descends, and consider a flat of type \(\finA[2]\) as before. Then we have a relation \(\flow[\sref[\ver]]\flow[\sref[\ver*]]\flow[\sref[\ver]]=\flow[\sref[\ver*]]\flow[\sref[\ver]]\flow[\sref[\ver*]]\), that is
  \begin{equation}\label{eqn:flat}
    \smut[\ver]<\pm>\smut[\ver*]<\pm>\smut[\ver]<\pm>=\smut[\ver*]<\pm>\smut[\ver]<\pm>\smut[\ver*]<\pm>.
  \end{equation}
  Now, by \cite{hirano_faithful_2018}, \cite{brav_braid_2011}, there is a monomorphism \(\inj{\Br[\finA[2]]}{\gAuteq{\Gcat*[\graph]}}\) (via \(\sref[\ver]\mapsto\smut[\ver]\)). Therefore \eqref{eqn:flat} implies that
  \begin{equation*}
    \sref[\ver]<\pm>\sref[\ver*]<\pm>\sref[\ver]<\pm>=\sref[\ver*]<\pm>\sref[\ver]<\pm>\sref[\ver*]<\pm>
  \end{equation*}
  holds in \(\Br[\finA[2]]\). \cref{lem:flow} says that this can only happen if there is a unique source and sink alcove around the flat. The \(\finA[1]\times\finA[1]\) case is similar.  
\end{proof}

\begin{lem}\label{lem:flowtriv}
  If \(\flow\) is a real flow then the monodromy monoid homomorphism
  \begin{equation*}
    \map{\End[\Path]{\alc}}{\Br\Gcat*[\graph]},
    \qquad \path\mapsto\flow[\path]
  \end{equation*}
  is trivial for all \(\alc\in\Alc\). 
\end{lem}
\begin{proof}
  We induct on the length of \(\map[\path]{\alc}{\alc}\). If this is two, then \(\path = (\alc\xrightarrow{\sref[\ver]}\alc*\xrightarrow{\sref[\ver]}\alc)\) for some \(\ver\in\Ver{\graph}\). By definition of a real functor we have \(\flow[\path] = \smut[\ver]<\pm>\smut[\ver]<\mp> = \id\). If \(\path\) has length greater than two then it must cross some hyperplane \(\hpln=\hpln[\weyl\sref[\ver]]\) twice. Thus we have a decomposition \(\path=\upgamma\circ\sref[\ver]\circ\upbeta\), where \(\map{\alc*}{C}\) is the first place \(\path\) crosses \(\hpln\) twice, and \(\upbeta\) does not cross any hyperplane twice. Then \(\upbeta\) is reduced.
  \begin{equation*}
    \begin{tikzpicture}
      \node (1) at (-1,-1) {$\alc$};
      \node (m-1) at (2,0.5) {$\alc*$};
      \node (m) at (2.5,-0.5) {$C$};
      \draw[dotted] (-2,0) -- (4,0);
      \node at (4.25,0) {$\hpln$};
      \draw[color=black,rounded corners=15pt,->] 
       (1)  --($(1)+(0.5,2)$) -- ($(1)+(1.5,1)$) --  node[left]{$\upbeta$}($(1)+(2,2.5)$)   -- (m-1);
      \draw[->] (m-1) to node[right] {$\sref[\ver]$} (m);
      \draw[->, bend left] (m) to node[below] {$\upgamma$} (1);
    \end{tikzpicture}
  \end{equation*}
  We may choose a reduced path \(\map[\upbeta']{\alc}{C}\), which does not cross \(\hpln\) because \(\alc\) and \(C\) lie on the same side of \(\hpln\), and any path between them crossing \(\hpln\) must do so twice and would not be reduced. Then the path \(\map{\alc}{\map{C}{\alc*}}\) is also reduced because it crosses each hyperplane at most once. Therefore \(\upbeta=\sref[\ver]\circ\upbeta'\) in the path category and so
  \begin{equation*}
    \flow[\path] = \flow[\upgamma]\flow[\map{\alc*}{C}]\flow[\upbeta] = \flow[\upgamma]\flow[\map{\alc*}{C}]\flow[\map{C}{\alc*}]\flow[\upbeta'] = \flow[\upgamma]\smut[\ver]<\pm>\smut[\ver]<\mp>\flow[\upbeta'] = \flow[\upgamma\circ\upbeta'],
  \end{equation*}
  which is of smaller length. The claim follows by induction.
\end{proof}

This explains the \say{real} in the terminology, with the monodromy morphism being trivial reflecting the lack of monodromy in a real hyperplane arrangement.

\begin{exm}\label{exm:abm}
  We give two examples of real flows on \(\Hpln[\graph]\) from \cref{exm:rkthreeaff} and \cref{exm:rkthreehyp}, the former has the most restrictions, while the latter has no restrictions.
  \begin{enumerate}
    \item If \(\graph = \affA[2]\), then we may recover the notion of direction on affine arrangements used in \cite{anno_stability_2015}, which involves choosing a connected component of the linearised arrangement. If we consider a flow \(\flow\) that has the source alcove (and hence the sink alcove) in the same position for every flat, then then this has the same effect as fixing the linear component in the same position. In Figure 3 we use the \say{top left} as the source around every \(\finA[2]\) flat, but the technology of flows allows us to construct many more examples by making compatible choices locally at flats.
    \item If \(\graph = \affA[2,3]\), then \(\Weyl[\graph]\) is the free group on three generators. This means that the associated Coxeter arrangement has no codimension two walls (as no root hyperplanes intersect inside the disc) and hence the condition on flats is vacuously satisfied, so any choice of real functor will automatically be a flow.
  \end{enumerate}
  \begin{figure}[ht]
    \begin{equation*}
      \begin{array}[c]{cc}
        \begin{overpic}[width=.48\textwidth]
          {image/A2affy.pdf}
          \put (10,41) {\rotatebox[origin=c]{90}{\huge$\textcolor{r}{\below}$}}
          \put (34,41) {\rotatebox[origin=c]{90}{\huge$\textcolor{r}{\below}$}}
          \put (60,41) {\rotatebox[origin=c]{90}{\huge$\textcolor{r}{\below}$}}
          \put (84,41) {\rotatebox[origin=c]{90}{\huge$\textcolor{r}{\below}$}}
          \put (22,63) {\rotatebox[origin=c]{90}{\huge$\textcolor{r}{\below}$}}
          \put (47,63) {\rotatebox[origin=c]{90}{\huge$\textcolor{r}{\below}$}}
          \put (73,63) {\rotatebox[origin=c]{90}{\huge$\textcolor{r}{\below}$}}
          \put (22,21) {\rotatebox[origin=c]{90}{\huge$\textcolor{r}{\below}$}}
          \put (47,21) {\rotatebox[origin=c]{90}{\huge$\textcolor{r}{\below}$}}
          \put (73,21) {\rotatebox[origin=c]{90}{\huge$\textcolor{r}{\below}$}}
          \put (40,74) {\rotatebox[origin=c]{150}{\huge$\textcolor{r}{\below}$}}
          \put (65,74) {\rotatebox[origin=c]{150}{\huge$\textcolor{r}{\below}$}}
          \put (27,52) {\rotatebox[origin=c]{150}{\huge$\textcolor{r}{\below}$}}
          \put (53,52) {\rotatebox[origin=c]{150}{\huge$\textcolor{r}{\below}$}}
          \put (77,52) {\rotatebox[origin=c]{150}{\huge$\textcolor{r}{\below}$}}
          \put (14,32) {\rotatebox[origin=c]{150}{\huge$\textcolor{r}{\below}$}}
          \put (40,32) {\rotatebox[origin=c]{150}{\huge$\textcolor{r}{\below}$}}
          \put (65,32) {\rotatebox[origin=c]{150}{\huge$\textcolor{r}{\below}$}}
          \put (27,10) {\rotatebox[origin=c]{150}{\huge$\textcolor{r}{\below}$}}
          \put (53,10) {\rotatebox[origin=c]{150}{\huge$\textcolor{r}{\below}$}}
          \put (27,74) {\rotatebox[origin=c]{210}{\huge$\textcolor{r}{\below}$}}
          \put (53,74) {\rotatebox[origin=c]{210}{\huge$\textcolor{r}{\below}$}}
          \put (14,52) {\rotatebox[origin=c]{210}{\huge$\textcolor{r}{\below}$}}
          \put (40,52) {\rotatebox[origin=c]{210}{\huge$\textcolor{r}{\below}$}}
          \put (66,52) {\rotatebox[origin=c]{210}{\huge$\textcolor{r}{\below}$}}
          \put (27,32) {\rotatebox[origin=c]{210}{\huge$\textcolor{r}{\below}$}}
          \put (53,32) {\rotatebox[origin=c]{210}{\huge$\textcolor{r}{\below}$}}
          \put (77,32) {\rotatebox[origin=c]{210}{\huge$\textcolor{r}{\below}$}}
          \put (40,10) {\rotatebox[origin=c]{210}{\huge$\textcolor{r}{\below}$}}
          \put (66,10) {\rotatebox[origin=c]{210}{\huge$\textcolor{r}{\below}$}}           
        \end{overpic} &
        \begin{overpic}[width=.43\textwidth,page=1]
          {image/A2hyp2J234.png}
          \put (32,50) {\rotatebox[origin=c]{150}{\huge$\textcolor{r}{\below}$}}
          \put (60,50) {\rotatebox[origin=c]{210}{\huge$\textcolor{r}{\below}$}}
          \put (46,35) {\rotatebox[origin=c]{270}{\huge$\textcolor{r}{\below}$}}
          \put (16,50) {\rotatebox[origin=c]{0}{\huge$\textcolor{r}{\below}$}}
          \put (74,50) {\rotatebox[origin=c]{0}{\huge$\textcolor{r}{\below}$}}
          \put (56,22) {\rotatebox[origin=c]{135}{\huge$\textcolor{r}{\below}$}}
          \put (32,22) {\rotatebox[origin=c]{225}{\huge$\textcolor{r}{\below}$}}
          \put (25,70) {\rotatebox[origin=c]{90}{\huge$\textcolor{r}{\below}$}}
          \put (65,70) {\rotatebox[origin=c]{270}{\huge$\textcolor{r}{\below}$}}
          \put (25,15) {\rotatebox[origin=c]{45}{\huge$\textcolor{r}{\below}$}}
          \put (37,7) {\rotatebox[origin=c]{90}{\huge$\textcolor{r}{\below}$}}
          \put (56,7) {\rotatebox[origin=c]{90}{\huge$\textcolor{r}{\below}$}}
          \put (63,15) {\rotatebox[origin=c]{-45}{\huge$\textcolor{r}{\below}$}}
          \put (80,45) {\rotatebox[origin=c]{150}{\huge$\textcolor{r}{\below}$}}
          \put (85,60) {\rotatebox[origin=c]{210}{\huge$\textcolor{r}{\below}$}}
          \put (75,75) {\rotatebox[origin=c]{30}{\huge$\textcolor{r}{\below}$}}
          \put (70,80) {\rotatebox[origin=c]{90}{\huge$\textcolor{r}{\below}$}}
          \put (30,80) {\rotatebox[origin=c]{90}{\huge$\textcolor{r}{\below}$}}
          \put (15,75) {\rotatebox[origin=c]{-30}{\huge$\textcolor{r}{\below}$}}
          \put (5,60) {\rotatebox[origin=c]{-30}{\huge$\textcolor{r}{\below}$}}
          \put (8,45) {\rotatebox[origin=c]{30}{\huge$\textcolor{r}{\below}$}}
        \end{overpic}
      \end{array}
    \end{equation*}
    \caption{Examples of real flows on \(\Hpln[\graph]\) for graphs \(\affA[2]\) and \(\affA[2,3]\).}
  \end{figure}
\end{exm}

  \section{Constructing real variations of stability}\label{sec:stab}
  Armed with a candidate for the real central charge from \cref{ssec:charge}, and a candidate for the alcove map from \cref{ssec:flow}, this section constructs real variations on \(\Gcat*[\graph]\), where \(\graph\) is a connected graph of affine or hyperbolic type without loops. From \eqref{tbl:hbolic}, the requirement of no loops is not very restrictive as it only excludes the families \(\Loop[\val]\), \(\Kronecker[\val, 1]\), \(\Kronecker[\val, 2]\), which do not admit interesting arrangements anyway.

Observe that when \(\graph\) has no loops, the Serre subcategory generated by \(\simp[\ver]\) coincides with \(\add[\simp[\ver]]\). This is because \(\Ext[\Ppa]{1}{\simp[\ver]}{\simp[\ver]} = 0\). To make the compatibility conditions of \cref{dfn:realstab}(2) easier to verify, we first simplify the description of our triangulated category and, in the process, answer an open question.

\subsection{Faithfulness of the standard heart}\label{ssec:faith}

The purpose of this subsection is to prove that
\begin{equation}\label{eqn:faith}
  \Db{\Nilp[\Ppa]} \iso \Db[\Nilp]{\Ppa} \defn* \Gcat*[\graph].
\end{equation}

\begin{dfn}\label{dfn:tfaith}
  Let \(\Abel\) be the heart of a bounded t-structure in a triangulated category \(\Gcat\). Then \(\Abel\) is called \emph{faithful} if the inclusion \(\Abel\subset\Gcat\) induces an equivalence of triangulated categories
  \begin{equation*}
    \map[\func]{\Db{\Abel}}{\Triang}.
  \end{equation*}
\end{dfn}

The functor \(\func\) above is often referred to as the \emph{realisation functor}. Since \eqref{eqn:faith} is equivalent to the standard heart \(\Nilp[\Ppa] \subset \Gcat*[\graph]\) being faithful, we need to show that
\begin{equation*}
  \map[\func]{\Db{\Nilp[\Ppa]}}{\Gcat*[\graph]}, \qquad \nilp\cx \mapsto \nilp\cx
\end{equation*}
is fully faithful and essentially surjective. The following result shows that, in our case, it is enough for \(\func\) to be fully faithful.

\begin{ppn}\label{ppn:esurj}
  Let \(\Abel\) be an abelian category, let \(\Serre\subseteq\Abel\) be a Serre subcategory, and suppose that the realisation functor
  \begin{equation*}
    \map[\func]{\Db{\Serre}}{\Db[\Serre]{\Abel}}     
  \end{equation*}
  is fully faithful. Then \(\func\) is essentially surjective.
\end{ppn}
\begin{proof}
  Since we are working in the bounded setting, we proceed by induction on the number of nonzero cohomology groups of \(\obj\in\Db[\Serre]{\Abel}\). If this is one, then \(\obj\simeq\coH{\deg}{\obj}\) for some \(\deg\in\Int\), which is in \(\Serre\) by definition. Hence \(\obj\in\image{\func}\). Now suppose that \(\coH{\sdot}{\obj} = 0\) in degrees less than some \(\deg\in\Int\) but \(\coH{\deg}{\obj} \neq 0\). Then using the \say{good} truncation functors \(\trunc{\deg}, \trunc*{\deg+1}\), there exists a triangle
  \begin{equation}\label{eqn:trunc}
    \begin{tikzcd}
      \trunc{\deg}\obj \ar[r] & \obj \ar[r] & \trunc*{\deg+1}\obj \ar[r] & \hspace{0em}
    \end{tikzcd}
  \end{equation}
  in \(\Db{\Ppa}\), see, for example, \cite[\S 3.4]{milicic_lectures_nodate}. As the t-structure restricts (\cref{lem:nilpres}), we have \(\trunc{\deg}\obj, \trunc*{\deg+1}\obj\in\Db[\Serre]{\Abel}\), which both have nonzero cohomology in fewer degrees than \(\obj\). Hence by induction,
  \begin{equation*}
    \trunc{\deg}\obj, \trunc*{\deg+1}\obj\in\image{\func}.
  \end{equation*}
  Moreover, as \(\func\) is fully faithful, \(\image{\func}\) is a triangulated subcategory of \(\Db[\Serre]{\Abel}\), so \eqref{eqn:trunc} forces \(\obj\in\image{\func}\).
\end{proof}

Since \(\Gcat*[\graph]\) is a full subcategory of \(\Db{\Ppa}\), we know that \(\func\) is fully faithful if and only if \(\map{\Db{\Nilp{\Ppa}}}{\Db{\Ppa}}\) is fully faithful. The following \devissage result gives a criterion to check this on objects in the Serre subcategory.

\begin{lem}[{\cite[4.2.13]{krause_homological_2022}}]\label{lem:ffaith}
  Let \(\Abel\) be an abelian category and \(\Abel*\) a full exact subcategory. Then the functor \(\map[\func]{\Db{\Abel*}}{\Db{\Abel}}\) is fully faithful if and only if for all \(\mod*, \mod\in\Abel*\) and \(\deg\in\Int[0]\), there is a bijection
  \begin{equation*}
    \map[\func[\deg]]{\Ext[\Abel*]{\deg}{\mod*}{\mod}}{\Ext[\Abel]{\deg}{\mod*}{\mod}}.
  \end{equation*}
\end{lem}

Recall that if \(\Abel\) is an abelian category and \(\mod*,\mod\in\Abel\), then there is an identification of \(\Ext[\Abel]{\deg}{\mod*}{\mod}\) with \(\Hom[\Db{\Abel}]{\mod*}{\mod\shift{\deg}}\), where the latter consists of equivalence classes of \emph{roof diagrams}
\begin{equation*}
  \begin{tikzcd}[column sep = 1em]
    & \ext\cx \ar[dl, "\qis"', "\sim"'{sloped}] \ar[dr, "\mor"] \\
    \mod* && \mod\shift{\deg},
  \end{tikzcd}
\end{equation*}
with \(\ext\cx\) a bounded complex in \(\Abel\), \(\qis\) a quasi-isomorphism and \(\mor\) is a chain map up to homotopy. Two such extensions \((\ext\cx, \qis, \mor), (\ext*\cx, \qis*, \mor*)\) are equivalent over \(\Abel\) if there exists a common \say{overroof}
\begin{equation*}
  \begin{tikzcd}[column sep = 1.5em]
    && \obj \ar[dl, "\qis**"', "\sim"'{sloped}] \ar[dr, "\mor**"] \\
    & \ext\cx \ar[dl, "\qis"', "\sim"'{sloped}] \ar[drrr, "\mor"'] && \ext*\cx \ar[dlll, "\qis*", "\sim"{sloped}] \ar[dr, "\mor*"] \\
    \mod* &&&& \mod\shift{\deg},
  \end{tikzcd}
\end{equation*}
where \(\qis**\) is a quasi-isomorphism and both squares commute in \(\Kb{\Abel}\). To distinguish between extensions over \(\Nilp[\Ppa]\) and extensions over \(\Mod{\Ppa}\), the latter having far more elements in a given equivalence class, we use the notation \(\brac{\ext\cx, \qis, \mor}[\Abel]\) for the equivalence class viewed over the category \(\Abel\).

We leverage the \(\Natural\)-grading of \(\Ppa\) and the resolution \eqref{eqn:projres} to exchange extensions over \(\Nilp[\Ppa]\) and \(\Mod{\Ppa}\). By \cite[10.4.7]{weibel_introduction_1994}, if \(\proj\cx\) is a bounded above complex of projectives in an abelian category \(\Abel\), then for all \(\obj\cx\in\Db*{\Abel}\), we have
\begin{equation}\label{eqn:weibel}
  \Hom[\Db*{\Abel}]{\proj\cx}{\obj\cx} \cong \Hom[\Kb*{\Abel}]{\proj\cx}{\obj\cx}.
\end{equation}

\begin{thm}\label{thm:faith}
  The standard t-structure with heart \(\Nilp[\Ppa]\) in \(\Gcat*[\graph]\) is faithful. That is, there is a fully faithful and essentially surjective functor
  \begin{equation*}
    \map*[\func]{\Db{\Nilp[\Ppa]}}{\Db[\Nilp]{\Ppa}}.
  \end{equation*}
\end{thm}
\begin{proof}
  In light of \cref{lem:ffaith}, we need to show that the map
  \begin{equation*}
    \map[\func[\deg]]{\Ext[\Nilp]{\deg}{\mod*}{\mod}}{\Ext[\Ppa]{\deg}{\mod*}{\mod}}, \qquad \brac{\ext\cx, \qis, \mor}[\Nilp] \mapsto \brac{\ext\cx, \qis, \mor}[\Ppa]
  \end{equation*}
  is a bijection for all \(\mod*,\mod\in\Nilp{\Ppa}\). By strong induction on module length, we only need to verify this for vertex simples \(\simp[\ver]\), \(\simp[\ver*]\) with \(\ver,\ver*\in\Ver{\graph}\). Indeed, for \(\mod*,\mod\in\Nilp[\Ppa]\), we can apply \(\Ext[\Nilp]{\deg}{\blank}{\mod}\) and \(\Ext[\Ppa]{\deg}{\blank}{\mod}\) to the short exact sequence
  \begin{equation*}
    \begin{tikzcd}
      0 \ar[r] & \mod*' \ar[r] & \mod* \ar[r] & \mod*'' \ar[r] & 0,
    \end{tikzcd}
  \end{equation*}
  where \(\mod*'\) has smaller length than \(\mod*\) and \(\mod*''\) is a vertex simple. This induces the long exact sequence
  \begin{equation*}
    \begin{tikzcd}[column sep = 1em]
      \cdots \ar[r] & \Ext[\Nilp]{\deg}{\mod*''}{\mod} \ar[r] \ar[d] & \Ext[\Nilp]{\deg}{\mod*}{\mod} \ar[r] \ar[d] & \Ext[\Nilp]{\deg}{\mod*'}{\mod} \ar[r] \ar[d] & \Ext[\Nilp]{\deg+1}{\mod*''}{\mod} \ar[r] \ar[d] & \cdots \\
      \cdots \ar[r] & \Ext[\Ppa]{\deg}{\mod*''}{\mod} \ar[r] & \Ext[\Ppa]{\deg}{\mod*}{\mod} \ar[r] & \Ext[\Ppa]{\deg}{\mod*'}{\mod} \ar[r] & \Ext[\Ppa]{\deg+1}{\mod*''}{\mod} \ar[r] & \cdots
    \end{tikzcd}
  \end{equation*}
  and so the five lemma applies since vertical maps are bijections for \(\mod*'\) and \(\mod*''\). Similarly, we can use a filtration of \(\mod\) and apply \(\Ext{\deg}{\mod*}{\blank}\).
  
  (Surjectivity). The Hom case \(\deg = 0\) is immediate, the case \(\deg = 1\) follows since \(\Nilp[\Ppa]\subset\Mod{\Ppa}\) is a Serre subcategory, and cases \(\deg > 2\) are trivial since \(\Db{\Ppa}\) is \(2\)-Calabi--Yau. Hence the only interesting case is \(\deg = 2\). Note that if \(\ver\neq\ver*\) then \(\Ext[\Ppa]{2}{\simp[\ver]}{\simp[\ver*]} = 0\) (also by the Calabi--Yau property) and there is nothing to prove, so consider \(\ver=\ver*\). Then \(\dim\Ext[\Ppa]{2}{\simp[\ver]}{\simp[\ver]} = \dim\Hom[\Ppa]{\simp[\ver]}{\simp[\ver]} = 1\) implies that any nonzero element of \(\Ext[\Ppa]{2}{\simp[\ver]}{\simp[\ver]} = \Hom[\Db{\Ppa}]{\simp[\ver]}{\simp[\ver]\shift{2}}\) is equivalent to a roof of the form
  \begin{equation*}
    \begin{tikzcd}[column sep = 1em]
      & \proj\cx \ar[dl, "\qis"', "\sim"'{sloped}] \ar[dr, "{\mor[\ver]}"] \\
      \simp[\ver] && \simp[\ver]\shift{2}
    \end{tikzcd}
  \end{equation*}
  with \(\proj\cx\) the projective resolution \eqref{eqn:projres}. To ease notation, we write \(\proj\) for \(\proj\cx[0] = \proj\cx[-2] = \proj[\ver]\), \(\proj*\) for \(\proj\cx[-1] = \proj*[\ver] = \direct*\proj[\ver*]\), and \(\simp\) for \(\simp[\ver]\). By noting that the differential in \(\proj\cx\) increases path length by one, this resolution can be made \(\Natural\)-graded. Indeed, the first three graded pieces are
  \begin{equation}\label{eqn:grproj}
    \begin{tikzcd}
      0 \arrow[r] & 0 \arrow[r] & 0 \arrow[r] & \grade{\proj}[0] \arrow[r] & \simp \arrow[r] & 0, \\
      0 \arrow[r] & 0 \arrow[r] & \grade{\proj*}[0] \arrow[r] & \grade{\proj}[1] \arrow[r] & 0 \arrow[r] & 0, \\
      0 \arrow[r] & \grade{\proj}[0] \arrow[r] & \grade{\proj*}[1]\arrow[r] & \grade{\proj}[2] \arrow[r] & 0 \arrow[r] & 0,
    \end{tikzcd}
  \end{equation}
  where for \(\mod\in\Mod{\Ppa}\), \(\grade{\mod}[\len]\) and \(\grade*{\mod}[\len]\) denote the \(\len\)-th graded piece and the direct sum of pieces in grades at least \(\len\), respectively. Since the differential on \(\proj\cx\) satisfies \(\diff[-2](\proj\rad)\subseteq\proj*\rad<2>\) and \(\diff[-1](\proj*\rad<2>)\inside\proj\rad<3>\) (this can be checked directly as \(\idemp[\ver]x\in\idemp[\ver]\rad=\proj\rad\) maps to \(\edge\idemp[\ver]x=\idemp[\ver*]\edge x\in\idemp[\ver*]\rad<2>\inside\proj*\rad<2>\), and is similarly for the other differential), \(\proj\cx\) descends to 
  \begin{equation}\label{eqn:approx}
    \begin{tikzcd}
      0 \ar[r] & \quo{\proj}{\proj\rad} \ar[r] \ar[d, equals] & \quo{\proj*}{\proj*\rad<2>} \ar[r] \ar[d, equals] & \quo{\proj}{\proj\rad<3>} \ar[r] \ar[d, equals] & \quo{\simp}{\simp\rad<4>} \ar[r] \ar[d, equals] & 0 \\
      0 \ar[r] & \quo{\proj}{\grade*{\proj}[1]} \ar[r] & \quo{\proj*}{\grade*{\proj*}[2]} \ar[r] & \quo{\proj}{\grade*{\proj}[3]} \ar[r] & \simp \ar[r] & 0.
    \end{tikzcd}
  \end{equation}
  Here, \(\proj\rad<\len> = \grade*{\proj}[\len]\) because \(\rad<\len>=\grade*{\Ppa}[\len]\) and \(\proj\) is a direct summand of \(\Ppa\). Thus, we have defined a \say{nilpotent approximation} of \(\proj\cx\) by
  \begin{equation*}
    \begin{tikzcd}
      \nilp\cx\defn (\quo{\proj}{\proj\rad} \ar[r] & \quo{\proj*}{\proj*\rad<2>} \ar[r] & \quo{\proj}{\proj\rad<3>}).
    \end{tikzcd}
  \end{equation*}
  This is clearly nilpotent in each degree. It is quasi-isomorphic to \(\simp\) concentrated in degree zero because it is equal to the direct sum of the three exact sequences in \eqref{eqn:grproj}. Any map \(\map{\proj\cx}{\simp}\) factors through \(\nilp\cx\) because \(\simp\rad<2>=0\). Therefore, \(\mor\) factors through \(\nilp\cx\) to give a map \(\mor'\) making
  \begin{equation*}
    \begin{tikzcd}[column sep = 1.5em]
      && \proj\cx \ar[dl, equals] \ar[dr, two heads] \\
      & \proj\cx \ar[dl, "\qis"', "\sim"'{sloped}] \ar[drrr, "\mor"'] && \nilp\cx \ar[dlll, "\qis'", "\sim"{sloped}] \ar[dr, "\mor'"] \\
      \simp &&&& \simp\shift{2}
    \end{tikzcd}
  \end{equation*}
  commute in \(\Kb{\Ppa}\). In fact, this diagram is already commutative in \(\Cb{\Ppa}\). Therefore
  \begin{equation*}
    \brac{\proj\cx, \qis, \mor}[\Ppa] = \brac{\nilp\cx, \qis', \mor'}[\Ppa] = \func[\deg]\brac{\nilp\cx, \qis', \mor'}[\Nilp].
  \end{equation*}

  (Injectivity). Now let \(\brac{\ext\cx, \qis, \mor}[\Nilp]\in\Ext[\Nilp]{\deg}{\simp[\ver]}{\simp[\ver*]}\) satisfy \(\brac{\ext\cx, \qis, \mor}[\Ppa] = 0\). By definition of the zero morphism in the derived category (see, for example, \cite[Lemma 2.1.4]{milicic_lectures_nodate}), there exists some \(\proj*\cx\in\Db{\Ppa}\) and a quasi-isomorphism \(\map*[\qis*]{\proj*\cx}{\ext\cx}\) such that \(\mor\circ\qis* = 0\). We may then replace \(\proj*\cx\) with the projective resolution \(\proj\cx\) of \(\simp[\ver]\) using \cite[Lemma 0649]{stacks_project}, which says that there exists a map of complexes \(\map[\qis**]{\proj\cx}{\proj*\cx}\), making the diagram
  \begin{equation*}
    \begin{tikzcd}
      \simp[\ver] & \proj*\cx \arrow[l, "{\qis\circ\qis*}"'] \\
      \proj\cx \arrow[u] \arrow[ur, dashed, "\qis**"']
    \end{tikzcd}
  \end{equation*}
  commute in \(\Cb{\Ppa}\), provided that \(\qis\circ\qis*\) is surjective in every degree. This is the case for us because \(\simp[\ver]\) is one-dimensional. Then \(\qis**\) is also a quasi-isomorphism. Therefore, there is a quasi-isomorphism \(\map[\qis*\circ\qis**]{\proj\cx}{\ext\cx}\) with \(\mor\circ(\qis*\circ\qis**) = 0\). Writing \(\mor*\defn\qis\circ\qis*\circ\qis**\), we have
  \begin{equation*}
    \brac{\ext\cx, \qis, \mor}[\Ppa] = \brac{\proj*\cx, \qis\circ\qis*, 0}[\Ppa] = \brac{\proj\cx, \mor*, 0}[\Ppa].
  \end{equation*}
  As in the surjectivity argument, we wish to realise the equivalence class \(\brac{\proj\cx, \mor*, 0}[\Ppa]\) as an element of \(\Ext[\Nilp]{\deg}{\simp[\ver]}{\simp[\ver*]}\) via a nilpotent approximation. Since \(\ext\cx\) is a bounded complex of nilpotent modules, there exists an \(\len\in\Int[0]\) such that \(\ext\cx[\deg]\rad<\len> = 0\) for all \(\deg\in\Int\). Taking \(\len*\) to be the greater of \(2\) and this \(\len\), we define the complex of nilpotent modules
  \begin{equation*}
    \begin{tikzcd}
      \nilp\cx\defn (\quo{\proj}{\proj\rad<\len*>} \ar[r] & \quo{\proj*}{\proj*\rad<\len*+1>} \ar[r] & \quo{\proj}{\proj\rad<\len*+2>}).
    \end{tikzcd}
  \end{equation*}
  Then \(\map[\qis*\circ\qis**]{\proj\cx}{\ext\cx}\) factors as \(\surj{\proj\cx}{\map{\nilp\cx}{\ext\cx}}\) in \(\Cb{\Ppa}\). Since the quotient map \(\surj[\prj]{\proj\cx}{\nilp\cx}\) and the composition \(\qis*\circ\qis**\) are both quasi-isomorphisms, the morphism \(\map*[\qis']{\nilp\cx}{\ext\cx}\) is a quasi-isomorphism. Thus the commutative diagram
  \begin{equation*}
    \begin{tikzcd}
      & & & \proj\cx \ar[dl, "\qis**"', "\sim"'{sloped}] \ar[dr, "\prj", two heads, shift left] \\
      & & \proj*\cx \ar[dl, "\qis*"', "\sim"'{sloped}] \ar[ddrrrr, "0"'] & & \nilp\cx \ar[dlll, "\qis'", "\sim"{sloped}] \ar[ddrr, "0"] \\
      & \ext\cx \ar[drrrrr, "\mor"'] \ar[dl, "\qis"', "\sim"'{sloped}] \\
      \simp[\ver] & & & & & & \simp[\ver*]\shift{\deg}
    \end{tikzcd}
  \end{equation*}
  implies that \((\mor\circ\qis')\circ\prj = 0\), as \(\mor\circ\qis*\circ\qis** = 0\). Now consider the chain map \(\map[\mor\circ\qis'\circ\prj]{\proj\cx}{\simp[\ver*]\shift{\deg}}\), which is null homotopic by the above. This killing homotopy clearly factors through \(\nilp\cx\) via
  \begin{equation*}
    \begin{tikzcd}[column sep = 3.5em]
      0 \ar[r] & \proj \ar[r] \ar[d, two heads] & \proj* \ar[r] \ar[d, two heads] & \proj \ar[r] \ar[d, two heads] & 0 &&[-2em] \proj\cx \ar[d, "\prj", two heads] \\
      0 \ar[r] & \quo{\proj}{\proj\rad<\len*>} \ar[r] \ar[d] \ar[dl, dashed, "\mor*'"] & \quo{\proj*}{\proj*\rad<\len*+1>} \ar[r] \ar[d] \ar[dl, dashed, "\mor*'"] & \quo{\proj}{\proj\rad<\len*+2>} \ar[r] \ar[d] \ar[dl, dashed, "\mor*'"] & 0 &&[-2em] \nilp\cx \ar[d, "\mor\circ\qis'"] \\
      0 \ar[r] \ar[from=uur, "\mor*"', crossing over, near start] & \simp[\ver]\shift{\deg} \ar[r] \ar[from=uur, "\mor*"', crossing over, near start] & 0 \ar[r] \ar[from=uur, "\mor*"', crossing over, near start] & 0 \ar[r] & 0 &&[-2em] \simp[\ver*]\shift{\deg}
    \end{tikzcd}
  \end{equation*}
  and \(\mor*'\) is a killing homotopy by the surjectivity of \(\prj\). Hence, \(\mor\circ\qis'=0\), and then there is a commutative diagram
  \begin{equation*}
    \begin{tikzcd}
      & & \nilp\cx \ar[dl, "\qis'"', "\sim"'{sloped}] \ar[dr, equals] \\
      & \ext\cx  \ar[dl, "\qis"', "\sim"'{sloped}] \ar[drrr, "\mor"'] & & \nilp\cx \ar[dlll, "\qis\circ\qis'", "\sim"{sloped}] \ar[dr, "0"] \\
      \simp[\ver] & & & & \simp[\ver*]\shift{\deg},
    \end{tikzcd}
  \end{equation*}
  proving that \(\brac{\ext\cx, \qis, \mor}[\Nilp] = \brac{\nilp\cx, \qis\circ\qis', 0}[\Nilp] = 0\).
\end{proof}

\begin{rmk}
  The above proof is specific to the preprojective algebra, but we expect this proof to hold in greater generality, provided the algebra admits a \(\Natural\)-grading and has finite global dimension.
\end{rmk}

As a consequence, all of the hearts we assign to alcoves are faithful.

\begin{cor}\label{cor:faith}
  If \(\func\) is an autoequivalence of \(\Gcat*[\graph]\), then \(\func(\Nilp[\Ppa])\) is a faithful heart.
\end{cor}
\begin{proof}
  By \cite[Theorem 3.7.3]{gelfand_homological_1994}, a necessary and sufficient condition for a heart \(\Abel\inside\Triang\) to be faithful is that \(\Hom[\Triang]{\obj}{\obj*\shift{\deg}}\) is generated by \(\Hom[\Triang]{\obj}{\obj*\shift{1}}\) for all \(\deg\in\Int[2]\) and all \(\obj,\obj*\in\Abel\). Since this condition must hold for \(\Nilp[\Ppa]\) by \cref{thm:faith}, applying \(\func\) gives the result for \(\func(\Nilp[\Ppa])\).
\end{proof}

The following general result for Serre subcategories is useful when considering \(\Abel[\hpln]\filt[1]\) and \(\Abel*[\hpln]\filt[1]\) in \cref{dfn:realstab}.

\begin{lem}\label{lem:prjmap}
  If \(\Serre\) is a Serre subcategory of an abelian category \(\Abel\), then there exists a fully faithful functor
  \begin{equation*}
    \map{\quo{\Abel}{\Serre}}{\quo{\Db{\Abel}}{\Db[\Serre]{\Abel}}}.
  \end{equation*}
\end{lem}
\begin{proof}
  First consider the exact sequence of abelian categories
  \begin{equation*}
    \begin{tikzcd}
      0 \ar[r] & \Serre \ar[r] & \Abel \ar[r, "\prj"] & \quo{\Abel}{\Serre} \ar[r] & 0,
    \end{tikzcd}
  \end{equation*}
  which, by \cite[Theorem 3.2]{miyachi_localization_1991}, induces an exact sequence of triangulated categories
  \begin{equation}\label{eqn:triseq}
    \begin{tikzcd}
      0 \ar[r] & \Db[\Serre]{\Abel} \ar[r] & \Db{\Abel} \ar[r] & \Db{\quo{\Abel}{\Serre}} \ar[r] & 0.
    \end{tikzcd}
  \end{equation}
  Now, the fully faithful functor \(\map{\quo{\Abel}{\Serre}}{\Db{\quo{\Abel}{\Serre}}}\) factors as
  \begin{equation*}
    \begin{tikzcd}
      \quo{\Abel}{\Serre} \ar[r] & \quo{\Db{\Abel}}{\Db[\Serre]{\Abel}} \ar[r, "\sim"] & \Db{\quo{\Abel}{\Serre}},
    \end{tikzcd}
  \end{equation*}
  where the second functor is an equivalence using \eqref{eqn:triseq}. Hence the first functor must be fully faithful.
\end{proof}

\begin{ppn}\label{ppn:gwyn}
  If \(\Abel\) is an abelian category and \(\Serre\inside\Abel\) is Serre, then \(\Db[\Serre]{\Abel} = \thick[\Serre]\) as thick triangulated subcategories of \(\Db{\Abel}\).
\end{ppn}
\begin{proof}
  One inclusion is immediate, by definition of \(\thick[\Serre]\). Now every \(\ext\cx\in\Db[\Serre]{\Abel}\) admits a filtration by its cohomology objects (see, for example, \cite[Proposition 3.3.1]{bayer_tour_2011}):
  \begin{equation*}
    \begin{tikzcd}[column sep = 2em]
      0 \ar[rr] && \obj\filt[a] \ar[dl] \ar[r] & \cdots\cdots \ar[r] & \obj\filt[b-1] \ar[rr] && \obj\filt[b] = \ext\cx \ar[dl] \\
      & \coH{-a}{\ext\cx}\shift{a} \ar[ul, dashed] &&&& \coH{-b}{\ext\cx}\shift{b} \ar[ul, dashed]
    \end{tikzcd}
  \end{equation*}
  where all triangles are exact and each \(\coH{\deg}{\ext\cx}\in\Serre\). Hence \(\ext\cx\) can be built out of \(\Serre\) in finitely many steps, that is, \(\ext\cx\in\thick[\Serre]\).
\end{proof}

\subsection{Preparation}\label{ssec:prep}

We recall the calibration of the previous sections, writing \(\fHrt\defn\Nilp[\Ppa]\) and \(\Gcat\defn\Gcat*[\graph]\) to ease notation.

\begin{stp}\label{stp:rstab}
  We fix a connected graph \(\graph\) of affine or hyperbolic type without loops. The Coxeter arrangement \(\Hpln\) consists of a hyperplane \(\hpln[\rt]\subset\Ths\) dual to each positive real root \(\rt = \weyl\rt[\ver]\), which we view intersected with the set
  \begin{equation*}
    \Level[\graph] = \set*{\stab\in\Ths}[\sum\nolimits_{\ver\in\Ver{\graph}}\stab[\ver] \geq 0,\ \stab\tran\paren{\Adj{\GCM[\graph]}}\stab = 1].
  \end{equation*}
  Each alcove \(\alc\subset\Level\minus\Hpln\) is of the form \(\weyl\fAlc\) for some \(\weyl\in\Weyl\), and each of the \(\rank\) hyperplanes bordering \(\alc\) are the reflecting hyperplanes \(\hpln[\weyl\rt[\ver]]\) of the reflections \(\weyl\sref[\ver]\weyl\inv\) for \(\ver\in\Ver{\graph}\). Hence a neighbouring alcove \(\alc*\) has the form
  \begin{equation}\label{eqn:refl}
    \alc* = \weyl\sref[\ver]\weyl\inv\alc = \weyl\sref[\ver]\fAlc.
  \end{equation}
  We wish to assign to \(\alc=\weyl\fAlc\) the heart \(\Abel\defn\flow[\weyl](\Nilp[\Ppa])\), where \(\flow\) is a choice of real flow on \(\Gcat*[\graph]\) and
  \begin{equation*}
    \flow[\weyl] = \smut[\ver[\len]]<\pm_{\len}>\cdots\smut[\ver[1]]<\pm_{1}> \in \Br\Gcat*[\graph],
  \end{equation*}
  where \(\sref[\ver[\len]]\cdots\sref[\ver[1]]\) is a reduced expression for \(\weyl\). The neighbouring alcove \(\alc*\) is then assigned the heart \(\Abel*\defn\flow[\weyl]\smut[\ver]<\pm>(\Nilp[\Ppa])\), the sign depending on whether \(\alc*\) is above \((-)\) or below \((+)\) \(\alc\) in with respect to \(\flow\).
\end{stp}

Since \(\Ths\) admits the contragredient \(\Weyl\)-action, \(\pair{\weyl\stab}{\weyl\Kthy{\mod}} = \pair{\stab}{\Kthy{\mod}}\) for all \(\weyl\in\Weyl\), \(\stab\in\Ths\), \(\mod\in\Gcat*[\graph]\). Moreover, repeatedly applying \eqref{eqn:wpair} yields
\begin{equation}\label{eqn:ktilt}
  \Kthy{\flow[\weyl](\mod)} = \weyl\Kthy{\mod}.
\end{equation}

The following results allow us to mostly consider the case where \(\alc\) is the standard alcove \(\fAlc\) and \(\alc*\) is adjacent across the hyperplane \(\hpln[\ver]\) for some \(\ver\in\Ver{\graph}\).

\begin{lem}\label{lem:filt}
  In the context of \cref{stp:rstab}, consider alcoves \(\alc=\weyl\fAlc\) and \(\alc*=\weyl\sref[\ver]\fAlc\) separated by the hyperplane \(\hpln = \hpln[\weyl\rt[\ver]]\). Let \(\Abel = \flow[\weyl](\fHrt)\) and \(\Abel* = \flow[\weyl]\smut[\ver]<\pm>(\fHrt)\) be the corresponding hearts. Then, with notation from \cref{dfn:realstab}, we have
  \begin{enumerate}
    \item \(\Abel*[\hpln]\filt[1] = \flow[\weyl]\smut[\ver]<\pm>\flow[\weyl]\inv(\Abel[\hpln]\filt[1])\),
    \item \(\Gcat*[\alc*,\hpln]\filt*[1]\defn\set{\obj\in\Gcat}[\coH[\Abel*]{\sdot}{\obj}\in\Abel*[\hpln]\filt[1]] = \flow[\weyl]\smut[\ver]<\pm>\flow[\weyl]\inv(\Gcat*[\alc,\hpln]\filt*[1])\),
    \item \(\Abel[\hpln]\filt[1] = \add[\flow[\weyl]\simp[\ver]]\) and \(\Abel*[\hpln]\filt[1] = \add[\flow[\weyl]\simp[\ver]\shift{\mp1}]\),
    \item \(\Gcat*[\alc,\hpln]\filt*[1] = \thick[\flow[\weyl]\simp[\ver]] = \thick[\flow[\weyl]\simp[\ver]\shift{\mp1}] = \Gcat*[\alc*,\hpln]\filt*[1]\),
    \item \(\surj[\prj[\alc,\hpln]]{\Gcat}{\quo{\Gcat}{\Gcat*[\alc,\hpln]\filt*[1]}}\) sends \(\Abel\) and \(\Abel*\) to \(\quo{\Abel}{\Abel[\hpln]\filt[1]}\) and \(\quo{\Abel*}{\Abel*[\hpln]\filt[1]}\), respectively.
  \end{enumerate}
\end{lem}
\begin{proof}
  For (1) and (2) we prove only one inclusion, as the opposite inclusions are analogous.
  \begin{enumerate}
    \item Suppose that \(\mod\in\Abel[\hpln]\filt[1]\), namely \(\mod\in\Abel\) and \(\pair{\stab}{\Kthy{\mod}} = 0\) for all \(\stab\in\hpln[\weyl\rt[\ver]]\). Then \(\mod*\defn\flow[\weyl]\smut[\ver]<\pm1>\flow[\weyl]\inv(\mod)\) is such that
    \begin{equation*}
      \pair{\stab}{\Kthy{\mod*}} = \pair{\stab}{\weyl\sref[\ver]w\inv\Kthy{\mod}} = \pair{\sref[\ver]\weyl\inv\stab}{\weyl\inv\Kthy{\mod}} = \pair{\weyl\inv\stab}{\weyl\inv\Kthy{\mod}} = 0,
    \end{equation*}
    since \(\sref[\ver]\weyl\inv\stab\in\hpln[\ver]\), and so \(\mod*\in\Abel*[\hpln]\filt[1]\).
    \item Suppose that \(\obj\in\Gcat*[\alc,\hpln]\filt*[1]\), so \(\coH[\Abel]{\deg}{\obj}\in\Abel[\hpln]\filt[1]\) for all \(\deg\in\Int\). Set \(\func\defn\flow[\weyl]\smut[\ver]<\pm>\flow[\weyl]\inv\), so that \(\func(\Abel) = \Abel*\). If we set \(\obj*\defn\func\obj\), then
    \begin{equation*}
      \coH[\Abel*]{\deg}{\obj*} \stackrel{\eqref{eqn:mutcoh}}{\iso} \func\coH[\Abel]{\deg}{\obj}\in\func(\Abel[\hpln]\filt[1])= \Abel*[\hpln]\filt[1],
    \end{equation*}
    where the final equality is just an application of (1). Thus \(\func(\Gcat*[\alc,\hpln]\filt*[1]) \subseteq \Gcat*[\alc*,\hpln]\filt*[1]\), as required.
    \item The Serre subcategory \(\Abel[\hpln]\filt[1]\) consists of objects \(\mod\in\Abel\) such that \(\pair{\blank}{\Kthy{\mod}}\) has a zero of order at least one along \(\hpln[\weyl\rt[\ver]]\). If \(\stab\in\hpln[\weyl\rt[\ver]]\) then \(\weyl\inv\stab\in\hpln[\ver]\). We also have \(\sref[\ver]\weyl\inv\stab\in\hpln[\ver]\) since \(\sref[\ver]\hpln[\ver] = \hpln[\ver]\). Writing \(\mod = \flow[\weyl]\nilp\) for some nilpotent module \(\nilp\), we have
    \begin{equation*}
      \pair{\stab}{\Kthy{\mod}} = \pair{\stab}{\Kthy{\flow[\weyl]\nilp}} = \pair{\stab}{\weyl\Kthy{\nilp}} = \pair{\weyl\inv\stab}{\Kthy{\nilp}},
    \end{equation*}
    which, if zero, forces \(\Kthy{\nilp}\in\Int\rt[\ver]\). Thus \(\Abel[\hpln]\filt[1] = \flow[\weyl](\add[\simp[\ver]]) = \add[\flow[\weyl]\simp[\ver]]\), and combining \cref{lem:filt}(1) with \cref{lem:smuts} gives the result for \(\Abel*\).
    \item The thick triangulated subcategory \(\Gcat*[\alc,\hpln]\filt*[1]\) consists of objects in \(\Gcat\) with cohomology in \(\Abel[\hpln]\filt[1] = \add[\flow[\weyl]\simp[\ver]]\). By \cref{ppn:gwyn}, this is just \(\thick[\flow[\weyl]\simp[\ver]]\), and so
    \begin{equation*}
      \Gcat*[\alc*,\hpln]\filt*[1] \stackrel{(2)}{=} \flow[\weyl]\smut[\ver]<\pm>\flow[\weyl]\inv(\Gcat*[\alc,\hpln]\filt*[1]) = \flow[\weyl]\smut[\ver]<\pm>\thick[\simp[\ver]] \stackrel{\ref{lem:smuts}}{=}  \flow[\weyl]\thick[\simp[\ver]\shift{\mp1}] = \thick[\flow[\weyl]\simp[\ver]] = \Gcat*[\alc,\hpln]\filt*[1],
    \end{equation*}
    as required.
    \item For the final claim, consider the following diagram where, for readability, we have suppressed \(\hpln\) in the notation.
    \begin{equation*} \hspace{100em minus 1fil}
      \begin{tikzcd}[column sep = .9em]
        & \thick[\flow[\weyl]\simp[\ver]] \ar[r, "\sim"] \ar[d, equal] & \thick[\flow[\weyl]\simp[\ver]] \ar[r, equal] \ar[d, equal] & \thick[\flow[\weyl]\simp[\ver]\shift{\mp1}] \ar[d, equal] & \thick[\flow[\weyl]\simp[\ver]\shift{\mp1}] \ar[d, equal] \ar[l, "\sim"'] \\
        & \Db[\Abel\filt[1]]{\Abel} \ar[d, hook] & \Gcat*[\alc]\filt[1] \ar[d, hook] & \Gcat*[\alc*]\filt[1] \ar[d, hook'] & \Db[\Abel*\filt[1]]{\Abel*} \ar[d, hook'] \\
        \Abel \ar[r, hook] \ar[d, mapsto] & \Db{\Abel} \ar[r, "\sim"] \ar[d, two heads] & \Gcat \ar[r, equal] \ar[d, two heads, "{\prj[\alc]}"'] & \Gcat \ar[d, two heads, "{\prj[\alc*]}"] & \Db{\Abel*} \ar[l, "\sim"'] \ar[d, two heads] & \Abel* \ar[l, hook'] \ar[d, mapsto] \\
        \quo{\Abel}{\Abel\filt[1]} \ar[r, hook, dotted] & \quo{\Db{\Abel}}{\Db[\Abel\filt[1]]{\Abel}} \ar[r, "\sim"] & \quo{\Gcat}{\Gcat*[\alc]\filt[1]} \ar[r, equal] & \quo{\Gcat}{\Gcat*[\alc*]\filt[1]} & \quo{\Db{\Abel*}}{\Db[\Abel*\filt[1]]{\Abel*}} \ar[l, "\sim"'] & \quo{\Abel*}{\Abel*\filt[1]} \ar[l, hook', dotted]
      \end{tikzcd} \hfilneg
    \end{equation*}
    By (4), \(\prj[\alc,\hpln] = \prj[\alc*,\hpln]\) since \(\Gcat*[\alc,\hpln]\filt*[1] = \Gcat*[\alc*,\hpln]\filt*[1]\). The dotted maps in the bottom left and bottom right are fully faithful by \cref{cor:faith} and \cref{lem:prjmap} applied to the Serre subcategories \(\Abel[\hpln]\filt[1]\) and \(\Abel*[\hpln]\filt[1]\), respectively. The result follows by commutativity.    
    \qedhere
  \end{enumerate}
\end{proof}

In particular, the above proof shows that \(\quo{\Abel}{\Abel[\hpln]\filt[1]}\) and \(\quo{\Abel*}{\Abel*[\hpln]\filt[1]}\) are abelian subcategories of \(\quo{\Gcat}{\Gcat*[\alc,\hpln]\filt*[1]}\), so it makes sense to compare them in \eqref{eqn:quoprop}.

\begin{lem}\label{lem:intfilt}
  For all \(\ver\in\Ver{\graph}\), we have
  \begin{equation*}
    \smut[\ver](\add\simp[\ver]) = \smut[\ver](\fHrt)\inter\Db[\add\simp[\ver]]{\Ppa}.
  \end{equation*}
\end{lem}
\begin{proof}
  Let \(\mod = \smut[\ver](\nilp)\) for some \(\nilp\in\fHrt\). If \(\nilp\in\add\simp[\ver]\), then \(\mod = \smut[\ver](\simp[\ver]<\direct\len*>) = \simp[\ver]<\direct\len*>\shift{1}\) has cohomology in \(\add\simp[\ver]\), so \(\mod\in\Db[\add\simp[\ver]]{\Ppa}\). This proves \(\subseteq\). Conversely, if \(\coH{\deg}{\mod}\in\add\simp[\ver]\) for all \(\deg\in\Int\), then we have
  \begin{equation*}
      \Kthy{\mod} = \Sum{\paren{-1}<\deg>}[\deg\in\Int]\Qdim{\coH{\deg}{\mod}} \in \Int\rt[\ver],
  \end{equation*}
  and hence \(\Kthy{\nilp} = \sref[\ver]<2>\Kthy{\nilp} = \sref[\ver]\Kthy{\mod} \in \Int\rt[\ver]\), which forces \(\nilp\in\add[\simp[\ver]]\).
\end{proof}

\subsection{Main result}\label{ssec:stab}

\begin{thm}\label{thm:rstab}
  Let \(\graph\) be a connected graph of affine or hyperbolic type and let \(\flow\) be a real flow on its Coxeter arrangement \(\Hpln[\graph]\). Then \((\charge, \hrt)\) is a real variation of stability on \(\Gcat*[\graph]\), where \(\charge\) is the inclusion of \(\Level\) into \(\Ths\) and \(\hrt\) assigns the heart \(\flow[\weyl](\Nilp[\Ppa])\) to the alcove \(\alc=\weyl\fAlc\in\Alc\).
\end{thm}
\begin{proof}
  The fact that \(\flow\) is a real flow means that the heart \(\Abel\defn\flow[\weyl](\fHrt)\) assigned to \(\alc\in\Alc\) is well-defined regardless of the reduced expression used for \(\weyl\). Indeed, by combining \cref{lem:flowtriv} and \cref{ppn:flow} one can assign to each \(\weyl\in\Weyl\) a functor \(\flow[\weyl] = \flow[\ver[\len]]\cdots\flow[\ver[1]]\). If \(\weyl=\sref[\ver*[\len]]\cdots\sref[\ver*[1]]\) is another reduced expression for \(\weyl\), then
  \begin{equation*}
    \sref[\ver[\len]]\cdots\sref[\ver[1]] = (\sref[\ver[\len]]\cdots\sref[\ver[1]]\sref[\ver*[1]]\inv\cdots\sref[\ver*[\len]]\inv)\sref[\ver*[\len]]\cdots\sref[\ver*[1]] = \upgamma\sref[\ver*[\len]]\cdots\sref[\ver*[1]]
  \end{equation*}
  for some \(\upgamma\in\End[\gDel]{\alc}\). Then by \cref{lem:flowtriv} the image of \(\upgamma\) under \(\flow\) is the trivial automorphism, so \(\flow[\weyl]\) is independent of the chosen expression. We verify the properties of \cref{dfn:realstab} in turn.
  \begin{enumerate}
    \item Let \(\mod\in\Abel\) be nonzero and let \(\stab\in\alc\). Then by construction there exists \(\stab*\in\fAlc\) and \(\nilp\in\fHrt\) such that \(\stab = \weyl\stab*\) and \(\mod = \flow[\weyl](\nilp)\).
    Now apply \eqref{eqn:ktilt} so that
    \begin{equation*}
      \charge[\stab]\Kthy{\mod} = \pair{\stab}{\Kthy{\mod}} = \pair{\weyl\stab*}{\Kthy{\flow[\weyl](\nilp)}} = \pair{\stab*}{\nilp}.
    \end{equation*}
    This is clearly positive as the standard alcove lives in the positive orthant of \(\Ths\) and \(\Kthy{\nilp}\) is based by simples with non-negative coefficients. Thus the real central charge \(\map[\charge]{\Level}{\Ths}\) satisfies the positivity property.
    \item Now suppose that \(\alc*\) is an alcove adjacent to \(\alc\) across the hyperplane \(\hpln[\rt]\). Following \cref{stp:rstab}, we have \(\rt = \weyl\rt[\ver]\) for some \(\ver\in\Ver{\graph}\), and the heart assigned to \(\alc*\) is then
    \begin{equation*}
      \Abel* \defn
      \begin{cases}
        \flow[\weyl]\smut[\ver](\fHrt) & \text{if } \alc \above \alc*, \\
        \flow[\weyl]\smut[\ver]\inv(\fHrt) & \text{if } \alc \below \alc*.
      \end{cases}
    \end{equation*}
    We need to prove the compatibility conditions of \cref{dfn:realstab} only for \(\deg\in\set{0,1}\), as the Serre subcategories \(\Abel[\hpln]\filt[\deg], \Abel*[\hpln]\filt[\deg]\) are zero for \(\deg > 1\). This is because \(\charge\) is the inclusion map and so the pairing is linear. Recall that conditions \eqref{eqn:intprop} and \eqref{eqn:quoprop} are, in order of increasing difficulty,
    \begin{equation*}
      \Abel*[\hpln]\filt[1] = \Abel[\hpln]\filt[1]\shift{\pm1}, \qquad \Abel*[\hpln]\filt[1] = \Abel*[\hpln]\inter\Gcat*[\alc,\hpln]\filt*[1], \qquad \quo{\Abel*}{\Abel*[\hpln]\filt[1]} = \quo{\Abel}{\Abel[\hpln]\filt[1]}. 
    \end{equation*}
    The first is \cref{lem:filt}(3), the second follows from \cref{lem:filt}(4) (or by applying \(\flow[\weyl]\) to \cref{lem:intfilt}), and by \cref{lem:filt}(5) the third condition reduces to checking
    \begin{equation}\label{eqn:prjmap}
      \prj[\alc,\hpln](\Abel) = \prj[\alc,\hpln](\Abel*).
    \end{equation} 
    To this end, consider the commutative diagram
    \begin{equation*}
      \begin{tikzcd}
        \Gcat \ar[r, "{\flow[\weyl]\inv}"] \arrow[d, "{\prj[\alc,\hpln]}"'] & \Gcat* \ar[d, "{\prj[\ver]}"] \\
        \quo{\Gcat}{\thick[\flow[\weyl]\simp[\ver]]} \ar[r, "{\cl{\flow[\weyl]\inv}}"'] & \quo{\Gcat}{\thick[\simp[\ver]]}
      \end{tikzcd}
    \end{equation*}
    of triangulated categories, where \(\prj[\ver]\defn\prj[\fAlc, \hpln[\ver]]\). Since the horizontal arrows are equivalences, it follows that \eqref{eqn:prjmap} follows from the equality \(\prj[\ver](\fHrt) = \prj[\ver](\smut[\ver]<\pm>\fHrt)\), because
    \begin{equation*}
      \cl{\flow[\weyl]\inv}(\prj[\alc,\hpln](\Abel)) = \prj[\ver](\flow[\weyl]\inv(\Abel)) = \prj[\ver](\fHrt),
    \end{equation*}
    and similarly for \(\Abel*\). This allows us to reduce to the case \(\weyl = 1\), for which we have
    \begin{align*}
      \alc &= \fAlc, & \alc* &= \sref[\ver]\fAlc, \\
      \Abel &= \fHrt, & \Abel* &= \smut[\ver]<\pm1>(\fHrt), \\
      \Abel[\hpln]\filt[1] &= \add[\simp[\ver]], & \Abel*[\hpln]\filt[1] &= \add[\simp[\ver]\shift{\pm1}],
    \end{align*}
    by \cref{lem:filt}(3).
    
    Let \(\mod\in\Abel*\). If we assume that \(\alc\below\alc*\) with respect to the flow \(\flow\), then \cref{lem:tiltmut} implies \(\mod\in\Tilt*[\ver]{\Abel}\) so that \(\Hom[\Abel]{\coH{0}{\mod}}{\simp[\ver]} = 0\), \(\coH{-1}{\mod}\in\add[\simp[\ver]]\), and \(\coH{\deg}{\mod}=0\) otherwise. Truncating as in \eqref{eqn:trunc}, we obtain the triangle
    \begin{equation*}
      \begin{tikzcd}[ar symbol/.style = {draw=none,"\textstyle#1" description,sloped}, iso/.style = {ar symbol={\cong}}]
        \trunc{-1}\mod \ar[r] & \mod \ar[r] & \trunc*{0}\mod \ar[r] & \phantom{.} \\
        \simp[\ver]<\direct\len*>\shift{1} \ar[u, iso] \ar[r] & \mod \ar[u, equals] \ar[r] & \nilp \ar[u, iso] \ar[r] & \phantom{.}
      \end{tikzcd}
    \end{equation*}
    for some \(\nilp\in\Tors[\ver]\subseteq\Nilp[\Ppa]\). Applying \(\prj[\ver]\) kills the first term in the quotient by \(\thick[\simp[\ver]]\). Therefore \(\prj[\ver](\mod)\) is isomorphic to \(\prj[\ver](\nilp)\), giving \(\prj[\ver](\Abel*)\subseteq\prj[\ver](\Abel)\). For the reverse inclusion, take \(\nilp\in\Abel\) and use left tilts with the triangle
    \begin{equation*}
      \begin{tikzcd}
        \trunc{0}\nilp \ar[r] & \nilp \ar[r] & \trunc*{1}\nilp \ar[r] & \phantom{.}
      \end{tikzcd}
    \end{equation*}
    The argument for both inclusions is analogous if instead \(\alc\above\alc*\). \qedhere
  \end{enumerate}
\end{proof}

  \bibliographystyle{alpha}
  \bibliography{RealStab}
\end{document}